\documentclass{article}%
\usepackage{amsmath}
\usepackage{amsfonts}
\usepackage{amssymb}
\usepackage{graphicx}
\usepackage{color}
\usepackage{hyperref}

\usepackage{algorithm}
\usepackage{algpseudocode}

\newtheorem{theorem}{Theorem}

\newtheorem{definition}[theorem]{Definition}

\newtheorem{lemma}[theorem]{Lemma}

\newtheorem{remark}[theorem]{Remark}

\newenvironment{proof}[1][Proof]{\noindent\textbf{#1.} }{\ \rule{0.5em}{0.5em}}

\newcommand{\black}{\color{black}}

\begin{document}

\title{Arnold Diffusion in the Full Three-Body Problem}

\author{Maciej J. Capi\'nski \smallskip\\Faculty of Applied Mathematics\\AGH University of Krak\'ow \\al. Mickiewicza 30, 30-059 Krak\'ow, Poland \bigskip\\Marian Gidea \smallskip\\Department of Mathematical Sciences\\Yeshiva University\\New York, NY 10016, USA
}

\maketitle

\begin{abstract}
The full three-body problem, on the motion of three celestial bodies under their mutual gravitational attraction, is one of the oldest unsolved problems in classical mechanics. The main difficulty comes from the presence of unstable and chaotic motions, which make long-term prediction impossible. In this paper, we show that the full three-body problem exhibits a strong form of instability known as Arnold diffusion.
We consider  the planar full three-body problem, formulated as a perturbation of both the Kepler problem and the planar circular restricted three-body problem. We show that the system exhibits Arnold diffusion, in the sense that there is a transfer of energy -- of an amount independent of the perturbation parameter -- between the Kepler problem and the restricted three-body problem.
Our argument  is based on the topological method of correctly aligned windows, which is implemented into a computer assisted proof.
We demonstrate that the approach can be applied to physically relevant masses of the bodies, choosing a Neptune-Triton-asteroid system as an example.
%This approach can be applied to physically relevant masses of the bodies, such as those in a Neptune-Triton-asteroid system.
In this case, we obtain explicit estimates for the range of the perturbation parameter and for the diffusion time.
\bigskip
 
\noindent\textbf{Mathematics Subject Classification (2020):}
37J25,  	%Stability problems for finite-dimensional Hamiltonian and Lagrangian systems
37J40,  	%Perturbations of finite-dimensional Hamiltonian systems, normal forms, small divisors, KAM theory, Arnol'd diffusion
65G40,  	%General methods in interval analysis
%65G20  	Algorithms with automatic result verification
%37N05  	%Dynamical systems in classical and celestial mechanics
70F07,  	%Three-body problems
70F15,  	%Celestial mechanics
70K44 	%Homoclinic and heteroclinic trajectories for nonlinear problems in mechanics
%34C37  	%Homoclinic and heteroclinic solutions to ordinary differential equations


\end{abstract}

%\tableofcontents

%TCIDATA{Version=5.00.0.2606}
%TCIDATA{LaTeXparent=0,0,MMFedit.tex}

\section{Introduction}
Since Poincar\'e, understanding stability and instability  in Hamiltonian systems has been
regarded as one of the most important problems in dynamical systems.
Poincar\'e  was particularly interested in how instability arises in the three-body problem, concerning the motion of three bodies interacting under newtonian gravity.
An important question is whether the effect of   instability  can  accumulate over time and lead to large effects. This question is relevant, for instance, in studying
the stability of the solar system.

%***Version A.***
In 1964 V.I. Arnold \cite{Arnold64} conjectured that instability is a pervasive phenomenon in `typical' Hamiltonian systems.
Arnold illustrated this idea with a special example of a Hamiltonian system of two-and-a-half-degrees of freedom, depending on two parameters.
When both parameters are set to zero,  the system consists of two rotators and is  fully integrable.
When the first parameter becomes nonzero, the system transforms into a rotator and a pendulum.
The energy of the rotator and that of the pendulum are conserved, and the system is still integrable (with a singularity). However,  the pendulum introduces hyperbolic dynamics into the system, i.e., there  is a pair of directions along which trajectories exhibit exponential expansion and contraction, respectively.
Then one considers the case when both parameters are non-zero.
Arnold showed  that when the second parameter is non-zero and  sufficiently small, the system exhibits `global' instability,
in the sense that there  exist
trajectories along which the energy of the rotator changes by an amount independent of the small parameter.
These trajectories are obtained along transition chains of invariant tori (consisting of heteroclinic connections between nearby tori).
The system also displays chaotic dynamics.
Arnold formulated the general conjecture as follows: \textit{I believe that
the mechanism of ``transition chains" which guarantees that nonstability in our
example is also applicable to the general case (for example, to the problem of
three bodies).}
Since then, the `Arnold diffusion problem' has been  one of the fundamental topics in Hamiltonian dynamics.

The  three-body problem shares some similarities with Arnold's example, but there are also significant differences.
Suppose  that the masses of the three bodies are $m_0,m_1,m_2$, with $m_0$  and $m_2$  viewed as parameters.
%When these parameters are set  to zero, i.e. $m_0=m_2=0$, these two bodies move under the gravity of the third body $m_1$
When $m_1>0$ and $m_0 = m_2 = 0$, the two bodies $m_0$ and $m_2$ move under the gravity of $m_1$ without influencing its motion,
so  the system consists of two Kepler problems, hence is fully integrable, as in Arnold's example.
%When we let one of the mass parameters be positive, say $m_0>0$, then the system consisting of $m_0$ and $m_1$ is a Kepler problem, and the motion of the infinitesimal mass $m_2$ under the  gravity of
 %$m_0$ and $m_1$ represents the restricted three-body problem.
%Both the energy of the Kepler problem and that of the  restricted three-body problem are conserved.
%\blue The energy of the Kepler problem is conserved.
%\red If the Keplerian orbit is circular, then the energy of the restricted three-body problem is also  conserved.
We then let one of the two mass parameters be positive, say $m_0 > 0$, and assume that $m_0$ and $m_1$ are in circular Keplerian motion.  Then the motion of the infinitesimal mass $m_2$
under the gravity of $m_0$ and $m_1$ is described by the circular restricted three-body problem.
The energy of the Kepler problem is conserved, and the energy of the circular restricted three-body problem is also conserved.
The  restricted three-body problem introduces  hyperbolic dynamics into the system, as well as chaotic dynamics.
In particular, it is not integrable.
When we let the second parameter be positive and small,  i.e., $m_{2}>0$,   we obtain the  full three-body
problem, whose equations of motion are given in \eqref{eqn:TBP}.
We will show that there exist large transfers of energy between the restricted three-body problem and the Kepler
problem.

Specifically,  we  consider a concrete model of the full three-body
problem, consisting of the  Neptune-Triton-asteroid system. The masses of
Neptune and Triton in normalized units are fixed and denoted  by  $ m_1= 1-\mu$ and $ m_0=\mu$, respectively, and the mass
of the asteroid is given by a small parameter  denoted $ m_2=\black\varepsilon ^{2}$.  The normalised mass of Triton is $\mu = 0.0002089$.
We assume that the initial conditions on Neptune and Triton  correspond to circular orbits around their common center of mass.
(If the asteroid was not present, the two bodies will maintain their circular orbits.)
In Section \ref{sec:Full-Kepler-R3BP} we show that, in an appropriate system of
coordinates and with the total energy  rescaled by $\varepsilon ^{2}$,   the full
three-body problem is described by  the Hamiltonian $\mathcal{\bar{H}}_{\varepsilon }$
\begin{equation}
\mathcal{\bar{H}}_{\varepsilon }=\bar{K}_{\varepsilon }+\bar{H}_{\varepsilon
}, \label{eq:H-as-2bp+rbp}
\end{equation}%
where $\bar{K}_{\varepsilon }$ is a perturbation of the Kepler problem and $%
\bar{H}_{\varepsilon }$ is a perturbation of the Planar Circular Restricted Three-Body
Problem (PCR3BP).
The expressions for $\mathcal{\bar{H}}_{\varepsilon }$, $\bar{K}_{\varepsilon }$, $\bar{H}_{\varepsilon
}$ will be provided in \eqref{eq:H-as-2bp+rbp-again}, \eqref{eq:K-bar-epsilon} and \eqref{eq:H-bar-epsilon}, respectively.
We note that both  $\bar{K}_{\varepsilon }$  and $\bar{H}_{\varepsilon }$ are quantities of order $O(1)$ with respect to the parameter $\varepsilon$.
For $\varepsilon =0$, both $\bar{K}_{0}$ and $\bar{H}_{0}$ are
constants of motion. For every $\varepsilon\ge 0$,  $\mathcal{\bar{H}}_{\varepsilon }$
is the total energy of the three-body problem, which is a constant of motion, but for $\varepsilon \neq 0$ the
quantities $%
\bar{K}_{\varepsilon }$ and $\bar{H}_{\varepsilon }$ no longer need to be
preserved along motions.

The main result of this paper is the following theorem, which we state here
in general terms, and reformulate precisely in Theorem \ref{th:main-exact} from Section \ref{sec:main-statement}.

\begin{theorem}[Main theorem]\label{th:main-intro}
Consider the planar full three-body problem, where the masses of the bodies are $m_0,m_1$ and $m_2=\varepsilon^2$. For
fixed values of $m_0$ and $m_1$, as specified below, and for a mass parameter  $\varepsilon ^{2}$, from arbitrarily close to
$0$ up to  some $\varepsilon_0^2$, there exist trajectories along which we have a  transfer
of energy between $\bar{K}_{\varepsilon }$ and $\bar{H}_{\varepsilon }$
in an amount that is independent of the size of $\varepsilon$.

As a model, we consider a Neptune–Triton–asteroid system, where the physical (not normalized) masses of Neptune and Triton are  $m_1=1.024 \times 10^{26}$ kg and $m_0=2.1389 \times 10^{22}$ kg, respectively, in which case we obtain a transfer of energy as above for \emph{all} asteroid physical masses $m_2$ from zero up to $10^6$ kilograms.
\end{theorem}

%The energy transfer  described in Theorem \ref{th:main-intro} is related to the  Arnold diffusion problem for Hamiltonian systems \cite{Arnold64}. This problem amounts to showing that  integrable Hamiltonian systems subjected to small perturbations of generic type have `diffusing orbits' along which the action variable changes by an amount independent of the smallness of the perturbation. Arnold illustrated this phenomenon for a pendulum-rotator system subject to small a perturbation  of  special type, and showed that there are diffusing orbits along which the action of the rotator (or, equivalently, the energy of the rotator) changes by an amount that is independent of the perturbation parameter. These orbits were obtained along transition chains of invariant tori (consisting of heteroclinic connections between nearby tori). Arnold claimed: ``I believe that this mechanism of instability is applicable to the general case (for example, to the problem of three bodies)".

Theorem \ref{th:main-intro} provides an answer to  Arnold's conjecture for the full three-body problem, for a system with real-life parameters. We note that in our case the unperturbed system is not fully  integrable, since we already assume that $\mu>0$. We show that one of the first integrals of the unperturbed system can undergo a change of order $O(1)$ under perturbation.

  For a system with realistic parameters, our result is beyond the reach of existing analytic or perturbative techniques.
Our approach combines topological methods with interval arithmetic validated numerics, yielding a rigorous, computer-assisted proof of our result.
For our computer assisted proof we have used the CAPD library \cite{MR4283203}. While we focus on a specific system here, our result serves as a proof of concept that the approach can be applied in settings close to real life. That is, we can consider different values for the masses of the two large bodies and a wider range of masses for the asteroid, at the expense of increased computational complexity.

Arnold's conjecture has received considerable interest and has seen significant progress over the years, including \cite{bolotin1999unbounded,DelshamsLS00,Treschev02c,Mather04,Treschev04, DelshamsLS06a,DelshamsLS06b,Piftankin2006,GelfreichT2008,DelshamsHuguet2009,ChengY09,Mather12,KaloshinZ15,bernard2016arnold,
cheng2019variational,
Treschev12,GideaL17,Gelfreich2017,gidea_marco_2017,kaloshin2020arnold,GideaLlaveSeara20-CPAM,GideaLlaveSeara20-DCDS}.

Arnold diffusion has been studied for various classes of systems \cite{ChierchiaG94}: (1)~A priori stable, in which the unperturbed system is fully integrable, and so  we have KAM stability, (2)~A priori unstable, in which the unperturbed system already has stable and unstable manifolds associated to a normally hyperbolic cylinder, and (3)~A priori chaotic, in which the unperturbed system has stable and unstable manifolds that intersect transversally, thus determining horseshoe dynamics. Our system is a priori chaotic and therefore does not exhibit KAM stability. In that case, one can use the intersections of the stable and unstable manifolds of the cylinder to prove Arnold diffusion. %(One does not need to construct diffusion between KAM tori as in the a priori stable case.)
Our approach is based on constructing certain sets positioned along the stable and unstable manifolds, however, our argument is topological and we do not need to explicitly establish these manifolds or their intersections.

Some of the existing results on Arnold diffusion are devoted to `generic systems' subjected to `arbitrarily small perturbations' of ‘generic type’, while other works
are focused on verifying diffusion in specific models.

Our work falls into the latter category, as we obtain a mechanism of diffusion that is applicable  to a concrete system (specifically, the Neptune-Triton-asteroid system), for an explicit perturbation, and for perturbation parameter values ranging from arbitrarily close to zero up to some explicit  cutoff  value. We note that this cutoff value is specific to our computer implementation, and diffusion may still occur beyond this value.

%Specifically, for the concrete parameters for the Neptune-Triton-asteroid full three-body problem.
%, we show that we obtain a change in the energy of the asteroid for a mass parameter  from arbitrarily close to zero up to a physically relevant value of $10^6$ kilograms.
%\marginpar{commented out more}

Another distinguished feature of our work is that we show Arnold diffusion in the full three-body problem, rather than in the restricted three-body problem, which has been considered in several papers, e.g.,  \cite{fejoz2016kirkwood,delshams2016arnold,CapinskiGL17,delshams2019global,clarke2022inner,guardia2023degenerate,clarke2024counterexample}. The model that we consider here is realistic, in the sense that the model parameters correspond to actual bodies in our solar system.
In contrast, other works concerning Arnold diffusion in  the full $n$-body problem  consider systems that satisfy certain scaling assumptions that may not be verified in  our solar system (e.g., the hierarchical regime, where the planets are increasingly separated, or the planetary regime, where the masses of the planets orbiting the Sun are arbitrarily small \cite{clarke2022inner}).

We choose to work with the Neptune-Triton-asteroid system because we establish the existence of  diffusing orbits for initial conditions so that orbits of the planet and its moon are close  to circular.
The Neptune-Triton system has the smallest eccentricity in the solar system, so this makes it a natural example to consider.
Triton is also of great interest to astronomers due to its retrograde orbit. It is hypothesised to be a captured Kuiper Belt object, and may have had a binary companion around 4 billion years ago.

While in this work we focus on this particular system,  our mechanism of diffusion can be applied to other systems.

Our work shares similar methods to those of   \cite{MR4544807}, which proves Arnold diffusion in the elliptic restricted  three-body problem for the  same system  (i.e., Neptune-Triton-asteroid). Specifically,  \cite{MR4544807} obtains orbits along which the energy  (measured in terms of the Hamiltonian of the restricted three-body problem) drifts by an amount independent of the smallness parameter, represented by the eccentricity of the elliptical orbits,
as well as orbits along which the energy makes chaotic excursions. It also provides  explicit estimates on the diffusion time and on the  Hausdorff dimension of the set of chaotic orbits.  Moreover, it shows that there are orbits
along which the time evolution of energy approaches a stopped  diffusion process (Brownian motion with drift), as  the perturbation parameter tends to $0$, and that any prescribed values of the drift and  variance for the limiting Brownian motion can be realized.
The methodology in \cite{MR4544807} relies on topological methods that are implemented in
computer assisted proofs.

The model considered in this paper presents several significant differences from \cite{MR4544807}. The smallness parameter $\varepsilon$  in this case is the square root of the mass of the smallest body (i.e., the asteroid). The diffusion is in terms of the  energy of the Hamiltonian of the restricted three-body problem $\bar{H}_{0}$  (see \eqref{eq:H-as-2bp+rbp}).
We show that there are orbits along which $\bar{H}_{0}$ drifts by $O(1)$, i.e., an amount independent of  $\varepsilon$. The corresponding Kepler energy  $\bar{K}_{0}$  along these orbits changes  by $O(\varepsilon^2)$.
Our construction of diffusing orbits starts with some initial conditions so that Neptune and Triton move on Kepler circular orbits about the center of mass, while the asteroid moves on a homoclinic  to a Lyapunov periodic orbit about the equilibrium point $L_1$ of the  restricted three-body problem.
Since the perturbation parameter is small, at any moment of time the motion of Neptune and Triton can be approximated by a `frozen' Kepler problem (whose orbits  become elliptic), and the motion of the asteroid by a `frozen' restricted three-body problem.
The effect of diffusion over a long period of  time is that the asteroid, which starts close to some Lyapunov periodic orbit, drifts along a homoclinic orbit to another Lyapunov periodic orbit, which is of different size (the size difference is $O(1)$), while the motion of the two massive bodies move along a Kepler orbit, whose semi-major axes changes slowly (the change is of order $O(\varepsilon)$).
%The effect of diffusion over a long period of  time is that the asteroid will move on a homoclinic orbit to some Lyapunov periodic orbit whose size differs from the original one by an large amount (of order $O(1)$),  while  the two massive bodies will move on Kepler orbits whose semi-major axes differ from the original ones  only by a small amount (of order $O(\varepsilon)$).

In our approach we take the following steps.

Starting with Newton's equations of motion we
first fix all classical integrals  and remove all
symmetries,   obtaining a reduced three-body problem. Then  we
use the methods of symplectic scaling and reduction as in \cite{MR1694376} to show that, in suitable coordinates,  the reduced
 three-body problem with one small mass is
given by a $3$-degree of freedom Hamiltonian that is the sum  of the Kepler problem and of the restricted three-body problem.
The Hamiltonian (representing the rescaled total energy of the three-body problem) is preserved, while the Kepler problem and the restricted three-body problem can transfer energy to one another.
The phase space of the system is $6$-dimensional and the solutions are restricted to a $5$-dimensional fixed energy level  of the three-body problem.

Second, we choose a  homoclinic orbit of the `frozen' restricted three-body problem  and construct a sequence of surfaces of section along it.
We  fix the total energy and   reduce the problem to a system  of section-to-section
$4$-dimensional maps.
We define a system of suitable coordinates on the surfaces of section, which we  denote by  $\left(u,s,\alpha,I\right)$, where $u$ will be an `unstable' coordinate, $s$  a `stable' coordinate, $\alpha $ an angle  coordinate, and $I$ an `action coordinate'  which measures the change in the Kepler energy starting from the initial level.
In the absence of the perturbation, for $\varepsilon=0$, the action $I$ is preserved and $\alpha $ is the same as  the time.
In the perturbed problem, for $\varepsilon>0$,  changes in  $I$ along  orbits of the section-to-section maps
 imply a transfer of energy between the Kepler problem and the restricted three-body problem.

Third, we use the  method of correctly aligned windows with cone conditions  as in \cite{MR4544807} to construct `connecting sequences' such that  the
action $I$ changes by $O(\varepsilon)$ along each connecting sequence. These connecting sequences have the property that the final window of a connecting sequence always overlaps with the initial window of another connecting sequence in such a way that there exist orbits of the section-to-section maps that pass through the
successive connecting sequences. By concatenating $O(1/\varepsilon)$ connecting sequences we can obtain  orbits along which the action $I$ changes by $O(1)$.  The  implementation of this method into a computer assisted proof yields the existence of diffusing orbits. It also provides  an explicit range for the values of the parameter $\varepsilon$,  as well as  explicit estimates for the diffusion time.

The paper is organised as follows. Section~\ref{sec:preliminaries} gives the setup for the full three-body problem.
Section~\ref {sec:Kepler-R3BP}  shows how to  reduce  the  full three-body problem to a  perturbation of the Kepler problem
and the circular restricted three-body problem. The statement of the main result of the paper, Theorem~\ref{th:main-exact} (which reformulates Theorem~\ref{th:main-intro}),  is given in Section~\ref{sec:main-statement}.
Section~\ref{sec:coordinates} introduces the aforementioned  surfaces of sections and the suitable local coordinates on them.
The method  of correctly aligned windows with cone conditions and connecting sequences is  reviewed in Section~\ref{sec:tools}.
The implementation of this method yields a  computer assisted proof  of Theorem~\ref{th:main-exact}, which is described in Section~\ref{sec:proof}.

%TCIDATA{Version=5.00.0.2606}
%TCIDATA{LaTeXparent=0,0,MMFedit.tex}

\section{Preliminaries -- the three-body problem}
\label{sec:preliminaries}
The problem describes the motion of three celestial bodies in space under
mutual Newtonian gravitational attraction. The equations of motions are given by
the Hamiltonian%
\begin{equation}\label{eqn:TBP}
H\left( q,p\right) =\sum_{i=0}^{2}\frac{\left\Vert p_{i}\right\Vert ^{2}}{%
2m_{i}}-\sum_{0\leq i<j\leq 2}\frac{\mathcal{G}m_{i}m_{j}}{\left\Vert
q_{i}-q_{j}\right\Vert },
\end{equation}%
where $m_{0},m_{1},m_{2}$ are the masses, $q_{0},q_{1},q_{2}$ are the
positions of the three bodies, $p_{0},p_{1},p_{2}$ are their momenta and $%
\mathcal{G}$ is the gravitational constant, with the symplectic form $dq\wedge dp$. The equations of motions are%
\begin{equation}
\dot{q}_{i}=\frac{\partial H}{\partial p_{i}},\qquad \dot{p}_{i}=-\frac{%
\partial H}{\partial q_{i}}\qquad \text{for }i=0,1,2.
\label{eq:full-3bp-ode-original}
\end{equation}%
By rescaling the masses we can consider that $\mathcal{G}=1$. 

Throughout the
paper we restrict our attention to the motion of the bodies on a plane, in
which case $q_{i},p_{i}\in \mathbb{R}^{2}$, for $i=0,1,2.$

When one of the masses is zero, say $m_{2}=0$, then the motion of the two
bodies with positive masses $m_{0}$ and $m_{1},$ which we refer to as
primaries, is described by the fully integrable Kepler problem. The motion
of the massless particle is described by the restricted three-body problem,
with the time dependent Hamiltonian
\begin{equation*}
\frac{\left\Vert p_{2}\right\Vert ^{2}}{2}-\frac{\mu }{r_{0}\left(
t,q_{2}\right) }-\frac{1-\mu }{r_{1}\left( t,q_{2}\right) },
\end{equation*}%
where $r_{0}\left( t,q_{2}\right) $ and $r_{1}\left( t,q_{2}\right) $ are
the distances of the massless particle to the first and the second primary,
respectively, and where $\mu =\frac{m_{0}}{m_{0}+m_{1}}$ is the normalised
mass of the first primary.

When $m_{2}=0$ and the primaries move on circular Keplerian orbits, by
passing to rotating coordinates the Hamiltonian of the restricted problem
becomes autonomous. This leads to the Planar Circular Restricted Three-Body
Problem (PCR3BP). We will show that in appropriate coordinates
%\footnote{%These coordinates are described in detail in section \ref{sec:Kepler-R3BP}.}
the Hamiltonian of the full three-body problem, for small $m_2>0$, can be expressed as a
perturbation of the sum of the Hamiltonians of the Kepler problem and of the
PCR3BP, with $m_{2}$ playing the role of the perturbation parameter. 

It is important to note that, when $m_2=0$, both the energy of  the Kepler problem and that of the PCR3BP  are
preserved. We will show that for $m_{2}>0$ we can
have a macroscopic transfer of energy between the Kepler problem and the
 PCR3BP, which is independent of the size of the perturbation.

\section{The full three-body problem as a perturbation of the Kepler problem
and the circular restricted three-body problem\label{sec:Kepler-R3BP}}

The aim of this section is to provide  an explicit formula for the planar full
three-body problem  in rotating coordinates, with appropriate
rescaling, so that it can be treated as a perturbation of
the Kepler problem and the Planar Restricted Three Body Problem (PR3BP). The
change of coordinates will also allow us to reduce the dimension of the
problem by making use of the classical conservation laws of the $n$-body
problem: conservation of the total linear momentum, the center of mass, and
the total angular momentum. These will allow us to reduce the dimension of
the system from the $12$-coordinates of the variables $\left( q,p\right) $
from (\ref{eq:full-3bp-ode-original}) down to  $6$-coordinates. This is done in
two steps.

In the first step, presented in section \ref{sec:rotating-Jacobi},  we derive
the Hamiltonian of the full three-body problem in the rotating Jacobi
coordinates. Our derivation combines the approach from the paper \cite%
{MR1694376} with a change to polar coordinates, which reduces the dimension from $12$ to $8$.
We perform
an additional step (\ref{eq:additional-change}), which is crucial for us,
but has not been considered in \cite{MR1694376}. Including (\ref%
{eq:additional-change}) allows us to reduce the dimension of the problem %by $2$,
from $8$ to $6$ by fixing the total angular momentum, while obtaining at
the same time explicit formulae for the equations of motions.

In the second step, presented in section \ref{sec:Full-Kepler-R3BP} we expand the coordinates in terms of the perturbation parameter
around a circular solution of the Kepler problem, which after
rescaling leads to the Hamiltonian (\ref{eq:H-as-2bp+rbp}). The expansion is
similar to the one from \cite{MR1694376}, but performed in our coordinates
from Section \ref{sec:rotating-Jacobi}.

\subsection{The Hamiltonian of the three--body problem in rotating Jacobi
coordinates\label{sec:rotating-Jacobi}}

\begin{figure}
	\begin{center}
		\includegraphics[height=3.5cm]{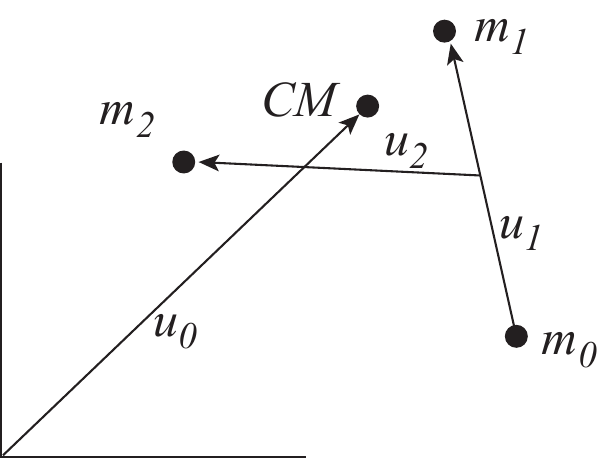}
	\end{center}
	\caption{Jacobi coordinates for the three body-problem.\label{fig:Jacobi}}
\end{figure}

We consider three particles with masses $m_{0},m_{1},m_{2}$, where $m_{2}$
will play the role of the smallest mass. We express the Hamiltonian of the planar three
body problem  by passing to the Jacobi coordinates \[(q_0,p_0,q_1,p_1,q_2,p_2)\mapsto (u_0,v_0,u_1,v_1,u_2,v_2).\]
The Jacobi coordinates are defined as follows \cite{MR1694376}:
\begin{equation}
\begin{array}{llll}
u_0=&\frac{m_0q_0+m_1q_1+m_2q_2}{m_0+m_1+m_2}, & v_0=&p_0+p_1+p_2, \\
u_1=&q_1-q_0,                                  & v_1=&-\frac{m_1}{m_0+m_1}p_0, \\
u_2=&q_2-\frac{m_0q_0+m_1q_1}{m_0+m_1},        & v_2=&-\frac{m_2}{m_0+m_1+m_2}(p_0+p_1+p_2)+p_2 .
\end{array}
\end{equation}
See Fig.~\ref{fig:Jacobi}.

By using  the conservation of the center of mass and of the total linear
momentum we can set $u_0=v_0=0$.
We thus obtain the Hamiltonian in Jacobi coordinates (see \cite{MR2468466}):
\begin{equation*}
H=\frac{\left\Vert v_{1}\right\Vert ^{2}}{2M_{1}}+\frac{\left\Vert
v_{2}\right\Vert ^{2}}{2M_{2}}-\frac{m_{0}m_{1}}{\left\Vert u_{1}\right\Vert
}-\frac{m_{1}m_{2}}{\left\Vert u_{2}-\alpha_{0}u_{1}\right\Vert }-\frac {%
m_{2}m_{0}}{\left\Vert u_{2}+\alpha_{1}u_{1}\right\Vert },
\end{equation*}
where the constants are%
\begin{align*}
M_{1} & =\frac{m_{0}m_{1}}{m_{0}+m_{1}},\qquad M_{2}=\frac{m_{2}\left(
m_{0}+m_{1}\right) }{m_{0}+m_{1}+m_{2}}, \\
\alpha_{0} & =\frac{m_{0}}{m_{0}+m_{1}},\qquad\alpha_{1}=\frac{m_{1}}{%
m_{0}+m_{1}}.
\end{align*}
The change to the Jacobi coordinates
reduces the dimension of the original problem (\ref{eq:full-3bp-ode-original}), from $(q,p)$ which is $12$-dimensional, to $(u,v)$ which is $8$-dimensional.

Our objective is to rewrite the equations so that, after letting $m_{2}\to 0$, we will obtain the Kepler problem in coordinates corresponding to $%
u_{1},v_{1}$, and the PCR3BP in coordinates corresponding to $u_{2},v_{2}$.

After passing to a set of coordinates that uniformly
rotate with frequency $\omega=1$ about the $z$-axis, the Hamiltonian takes form%
\begin{align*}
H & =\frac{\left\Vert v_{1}\right\Vert ^{2}}{2M_{1}}-u_{1}^{T}Jv_{1}-\frac{%
m_{0}m_{1}}{\left\Vert u_{1}\right\Vert } \\
& \quad+\frac{\left\Vert v_{2}\right\Vert ^{2}}{2M_{2}}-u_{2}^{T}Jv_{2}-%
\frac{m_{1}m_{2}}{\left\Vert u_{2}-\alpha_{0}u_{1}\right\Vert }-\frac {%
m_{2}m_{0}}{\left\Vert u_{2}+\alpha_{1}u_{1}\right\Vert },
\end{align*}
where
\begin{equation*}
J=\left(
\begin{array}{cc}
0 & 1 \\
-1 & 0%
\end{array}
\right),
\end{equation*}
with the symplectic form $du \wedge dv$.

The Hamiltonian $H$ is an integral of motion of the system%
\begin{equation*}
\dot{u}_{i}=\frac{\partial H}{\partial v_{i}},\qquad \dot{v}_{i} =-%
\frac{\partial H}{\partial u_{i}}\qquad\text{for }i=1,2.
\end{equation*}
Also the angular momentum, which is expressed as%
\begin{equation*}
A=u_{1}\times v_{1}+u_{2}\times v_{2},
\end{equation*}
is an integral of motion.

We now denote the masses of the primaries and  the small mass, respectively,  by
\begin{equation*}
m_{0}=\mu,\qquad m_{1}=1-\mu,\qquad m_{2}=\varepsilon^{2},
\end{equation*}
and let
\begin{equation*}
\nu=\mu\left( 1-\mu\right) .
\end{equation*}
In our convention $m_0<m_1$. 
The quantity $\varepsilon$ plays the role of the perturbation parameter. We introduce
the normalised masses $\mu$ and $1-\mu$ so that the equations that
follow  reduce to  the Kepler problem and the restricted
three-body problem when $m_2\to 0$. With this notation the constants become
\begin{align*}
M_{1} & =\nu,\qquad M_{2}=\frac{\varepsilon^{2}}{1+\varepsilon^{2}}, \\
\alpha_{0} & =\mu,\qquad\alpha_{1}=1-\mu,
\end{align*}
and the Hamiltonian takes the form%
\begin{equation*}
H=K+\tilde{H},
\end{equation*}
where
\begin{align*}
K & =\frac{\left\Vert v_{1}\right\Vert ^{2}}{2\nu}-u_{1}^{T}Jv_{1}-\frac {\nu%
}{\left\Vert u_{1}\right\Vert }, \\
\tilde{H} & =\left( \frac{1}{2\varepsilon^{2}}+\frac{1}{2}\right) \left\Vert
v_{2}\right\Vert ^{2}-u_{2}^{T}Jv_{2}-\frac{\varepsilon^{2}\left(
1-\mu\right) }{\left\Vert u_{2}-\mu u_{1}\right\Vert }-\frac{\varepsilon
^{2}\mu}{\left\Vert u_{2}+\left( 1-\mu\right) u_{1}\right\Vert },
\end{align*}
with the symplectic form $ du \wedge dv$.

After scalings
\[u_{i}\rightarrow u_{i},\,  v_{i}\rightarrow\nu v_{i},\,
K\rightarrow\nu^{-1}K,\,\tilde{H}\rightarrow\nu^{-1}\tilde{H},\,
\varepsilon^{2}\nu^{-1}\rightarrow\varepsilon^{2}\]
we can rewrite the above
as
\begin{align}
K & =\frac{1}{2}\left\Vert v_{1}\right\Vert ^{2}-u_{1}^{T}Jv_{1}-\frac {1}{%
\left\Vert u_{1}\right\Vert },  \label{eq:K-initial} \\
\tilde{H} & =\frac{1+\nu\varepsilon^{2}}{2\varepsilon^{2}}\left\Vert
v_{2}\right\Vert ^{2}-u_{2}^{T}Jv_{2}-\frac{\varepsilon^{2}\left(
1-\mu\right) }{\left\Vert u_{2}-\mu u_{1}\right\Vert }-\frac{\varepsilon
^{2}\mu}{\left\Vert u_{2}+\left( 1-\mu\right) u_{1}\right\Vert }.
\label{eq:H-initial}
\end{align}
We note that the $K$ is the Hamiltonian of the Kepler problem.
We also remark that $\tilde{H}$ is not defined for $\varepsilon=0$.

We now perform the symplectic change to polar coordinates
\[(u_1,v_1,u_2,v_2)\mapsto   (r_1, R_1,r_2,R_2, \theta_1,\Theta_1,\theta_2,\Theta_2)\]
given by:
\begin{align*}
u_{1} & =\left( r_{1}\cos\theta_{1},r_{1}\sin\theta_{1}\right) ,\qquad
u_{2}=\left( r_{2}\cos\theta_{2},r_{2}\sin\theta_{2}\right) , \\
v_{1} & =\left( R_{1}\cos\theta_{1}-\frac{\Theta_{1}}{r_{1}}\sin\theta
_{1},R_{1}\sin\theta_{1}+\frac{\Theta_{1}}{r_{1}}\cos\theta_{1}\right) , \\
v_{2} & =\left( R_{2}\cos\theta_{2}-\frac{\Theta_{2}}{r_{2}}\sin\theta
_{2},R_{2}\sin\theta_{2}+\frac{\Theta_{2}}{r_{2}}\cos\theta_{2}\right) ,
\end{align*}
which leads to%
\begin{align}
K & =\frac{1}{2}\left( R_{1}^{2}+\left( \frac{\Theta_{1}}{r_{1}}\right)
^{2}\right) -\Theta_{1}-\frac{1}{r_{1}},  \label{eq:K-rotating-Kepler} \\
\tilde{H} & =\frac{1+\nu\varepsilon^{2}}{2\varepsilon^{2}}\left(
R_{2}^{2}+\left( \frac{\Theta_{2}}{r_{2}}\right) ^{2}\right) -\Theta _{2}
\notag \\
& \quad-\frac{\varepsilon^{2}\left( 1-\mu\right) }{\sqrt{r_{2}^{2}+\mu
^{2}r_{1}^{2}-2\mu r_{1}r_{2}\cos\left( \theta_{2}-\theta_{1}\right) }}
\notag \\
& \quad-\frac{\varepsilon^{2}\mu}{\sqrt{r_{2}^{2}+\left( 1-\mu\right)
^{2}r_{1}^{2}+2\left( 1-\mu\right) r_{1}r_{2}\cos\left(
\theta_{2}-\theta_{1}\right) }},  \notag
\end{align}
with the symplectic form $dr \wedge dR + d \theta \wedge d \Theta$.

Note that the Hamiltonian only depends on the difference of the
polar angles $\theta_{2}-\theta_{1}$.

The angular momentum in the polar coordinates is
\[ A=u_1\times v_1+u_2\times v_2= \Theta_1+\Theta_2. \]

%We note that $K$ in (\ref{eq:K-rotating-Kepler}) is the Hamiltonian of the Kepler problem in the rotating polar coordinates.

We now   reduce the dimension of the system by fixing
the angular momentum. For this purpose we apply a symplectic change of
coordinates  (see \cite{celletti2010stability})
\[ (r_1, R_1,r_2,R_2,\theta_1,\Theta_1,\theta_2,\Theta_2)\mapsto (r_1,R_1, r_2,R_2,\phi_1,\Phi_1,\phi_2,\Phi_2)\]
given by
\begin{equation}
\begin{array}{lll}
\phi_{1}=\theta_{1}, & \qquad & \phi_{2}=\theta_{2}-\theta_{1},\medskip \\
\Phi_{1}=\Theta_{1}+\Theta_{2}, &  & \Phi_{2}=\Theta_{2},%
\end{array}
\label{eq:additional-change}
\end{equation}
which results in%
\begin{align}
K & =\frac{1}{2}\left( R_{1}^{2}+\left( \frac{\Phi_{1}-\Phi_{2}}{r_{1}}%
\right) ^{2}\right) -\left( \Phi_{1}-\Phi_{2}\right) -\frac{1}{r_{1}},
\notag \\
\tilde{H} & =\frac{1+\nu\varepsilon^{2}}{2\varepsilon^{2}}\left(
R_{2}^{2}+\left( \frac{\Phi_{2}}{r_{2}}\right) ^{2}\right) -\Phi_{2}
\label{eq:tilde-H} \\
& \quad-\frac{\varepsilon^{2}\left( 1-\mu\right) }{\sqrt{r_{2}^{2}+\mu
^{2}r_{1}^{2}-2\mu r_{2}r_{1}\cos\left( \phi_{2}\right) }}  \notag \\
& \quad-\frac{\varepsilon^{2}\mu}{\sqrt{r_{2}^{2}+\left( 1-\mu\right)
^{2}r_{1}^{2}+2\left( 1-\mu\right) r_{1}r_{2}\cos\left( \phi_{2}\right) }},
\notag
\end{align}
with the symplectic form $dr \wedge dR + d \phi \wedge d \Phi$.

The  total angular momentum of the three-body
problem  in the new coordinates is
\[ A= \Theta_1+\Theta_2=\Phi_1, \]
 which is a constant of motion. The coordinate $\phi_{1}$ is
ignorable; it does not influence the motion of the other variables. From now
on we restrict to $\Phi_{1}=1$, which leads to
\begin{equation}
K=\frac{1}{2}\left( R_{1}^{2}+\left( \frac{1-\Phi_{2}}{r_{1}}\right)
^{2}\right) -\left( 1-\Phi_{2}\right) -\frac{1}{r_{1}},   \label{eq:K-at-one}
\end{equation}
and ignore $\phi_{1}$.
By ignoring $\phi_1$ and setting $\Phi_1 = 1$ the system becomes $6$-dimensional, with the coordinates $(r_1, R_1,r_2,R_2,\phi_2,\Phi_2)$, and with the symplectic form $dr_1 \wedge dR_1 + dr_2 \wedge dR_2 + d\phi_2 \wedge d\Phi_2$. 
%This system depends on the coordinates $(r_1, R_1,r_2,R_2,\theta_2,\Theta_2)$ so is $6$-dimensional.

\subsection{The three-body problem as a perturbation of the Kepler and the circular
restricted three-body problem\label{sec:Full-Kepler-R3BP}}

We now perform a change of coordinates near initial conditions for $m_0$ and $m_1$  corresponding to the circular Kepler problem (that is, if $m_2$ was not present, starting with those initial conditions, $m_0$ and $m_1$ will move on circular orbits about their common center of mass). In the coordinates \eqref{eq:additional-change}, the circular Kepler problem corresponds to $r_1=1$ and $R_1=0$.
The new coordinate change (recalling that $\phi_1$ is ignored and $\Phi_1=1$) is
\[(r_1, R_1,r_2, R_2, \phi_2,\Phi_2)\mapsto (\bar{r}_1,\bar{R}_1,\bar{r}_2,\bar{R}_2, \bar{\phi}_2,\bar{\Phi}_2)\]
 given by:
\begin{equation}
\begin{array}{lllll}
r_{1}=1+\varepsilon \bar{r}_{1}, &  & r_{2}=\bar{r}_{2}, &  & \phi _{2}=\bar{%
\phi}_{2},\medskip \\
R_{1}=\varepsilon \bar{R}_{1}, &  & R_{2}=\varepsilon ^{2}\bar{R}_{2}, &  &
\Phi _{2}=\varepsilon ^{2}\bar{\Phi}_{2},%
\end{array}
\label{eq:eps-symplectic-scaling}
\end{equation}%
which is conformally symplectic with a conformal factor $\varepsilon ^{2}$. Our objective will
be to rewrite the Hamiltonian in the coordinates $\left( \bar{r}_{1},\bar{R}
_{1},\bar{r}_{2},\bar{R}_{2},\bar{\phi}_{2},\bar{\Phi}_{2}\right) $. The
coordinate change is motivated by the fact that for $\varepsilon =0$ and
$\Phi _{2}=0$, we have that $r_{1}=1,R_{1}=0$ describes  a circular solution of the Kepler
problem. So, the change of coordinates (\ref{eq:eps-symplectic-scaling}) can
be interpreted as an expansion of the coordinates in terms of the perturbation parameter around this circular solution.
We will show that (\ref{eq:eps-symplectic-scaling}) allows us to express the
Hamiltonian of the full three-body problem as a perturbation of the Kepler
problem and the PCR3BP.

We start with a lemma, which we prove in the Appendix~\ref%
{sec:proof-K-reduced}.

\begin{lemma}
\label{lem:K-reduced}The change of coordinates (\ref%
{eq:eps-symplectic-scaling}) leads to%
\begin{equation}
K=\varepsilon ^{2}\bar{K}_{\varepsilon }-\frac{3}{2},
\label{eq:K-in-bar-coordinates}
\end{equation}%
where%
\begin{equation}
\bar{K}_{\varepsilon }=\frac{1}{2\left( 1+\varepsilon \bar{r}_{1}\right) ^{2}%
}\left( \bar{R}_{1}^{2}+\bar{r}_{1}^{2}+2\varepsilon \bar{r}_{1}\left( 2\bar{%
\Phi}_{2}+\bar{R}_{1}^{2}\right) +\varepsilon ^{2}\left( \bar{\Phi}_{2}^{2}+%
\bar{r}_{1}^{2}\left( 2\bar{\Phi}_{2}+\bar{R}_{1}^{2}\right) \right) \right)
.  \label{eq:K-bar-epsilon}
\end{equation}
\end{lemma}

\begin{remark}
The interesting aspect of Lemma \ref{lem:K-reduced} is that (\ref%
{eq:K-in-bar-coordinates}) does not have any terms of order $\varepsilon $.
The details of the proof given in the appendix might obscure the simple
reason why this is so; which is that
\begin{align*}
\frac{1}{r_{1}}& =1+\varepsilon \bar{r}_{1}+O(\varepsilon ^{2}), \\
\left( \frac{1-\Phi _{2}}{r_{1}}\right) ^{2}& =1-2\varepsilon \bar{r}%
_{1}+O(\varepsilon ^{2}),
\end{align*}%
inserted into (\ref{eq:K-at-one}) clearly leads to the cancellation of the $%
\varepsilon $ terms.
\end{remark}

\begin{remark}
\label{rem:K0-rotor} Note that $\bar{K}_{0}$ is the harmonic oscillator
\begin{equation}
\bar{K}_{0}=\frac{1}{2}\left( \bar{R}_{1}^{2}+\bar{r}_{1}^{2}\right) .\label{eq:K_0-formula}
\end{equation}
\end{remark}

We can now formulate the following theorem.

\begin{theorem}
The equations of motion in coordinates $\left( \bar{r}_{1},\bar{R}_{1},\bar{r}_{2},\bar{R}_{2},\bar{\phi}_{2},\bar{\Phi}_{2}\right) $ are given by the
Hamiltonian%
\begin{equation}
H=-\frac{3}{2}+\varepsilon^{2}\left( \bar{K}_{\varepsilon}+\bar {H}%
_{\varepsilon}\right) ,
\label{eqn:H_final}
\end{equation}
where $\bar{K}_{\varepsilon}$ is given in (\ref{eq:K-bar-epsilon}) and
\begin{align}
\bar{H}_{\varepsilon} & =\frac{\left( 1+\nu\varepsilon^{2}\right) }{2}\left(
\bar{R}_{2}^{2}+\left( \frac{\bar{\Phi}_{2}}{\bar r_{2}}\right) ^{2}\right) -
\bar{\Phi}_{2} \notag \\
& \quad-\frac{1-\mu}{\sqrt{\bar{r}_{2}^{2}+\mu^{2}\left( 1+\varepsilon \bar{r%
}_{1}\right) ^{2}-2\mu\bar{r}_{2}\left( 1+\varepsilon\bar{r}_{1}\right)
\cos\left( \phi_{2}\right) }}  \notag \\
& \quad-\frac{\mu}{\sqrt{\bar r_{2}^{2}+\left( 1-\mu\right) ^{2}\left(
1+\varepsilon\bar{r}_{1}\right) ^{2}+2\left( 1-\mu\right) \bar{r}_{2}\left(
1+\varepsilon\bar{r}_{1}\right) \cos\left( \phi_{2}\right) }},  \label{eq:H-bar-epsilon}
\end{align}
with the symplectic form $d\bar r_1 \wedge d \bar R_1 + d \bar r_2 \wedge d \bar R_2 + d \bar \phi_2 \wedge d \bar \Phi_2$.
\end{theorem}

\begin{proof}
The result follows directly from Lemma \ref{lem:K-reduced} and by using (\ref%
{eq:eps-symplectic-scaling}) in (\ref{eq:tilde-H}).
\end{proof}

After dropping the constant term in  the Hamiltonian \eqref{eqn:H_final}  and then rescaling it by $
\varepsilon^{2} $ we obtain the Hamiltonian
\begin{equation}
\bar{\mathcal{H}}_{\varepsilon}=\varepsilon^{-2} (H + 3/2)=\bar{K}_{\varepsilon}+\bar{H}_{\varepsilon},
\label{eq:H-as-2bp+rbp-again}
\end{equation}
which is the main focus of our investigation.

\begin{remark}
The term $\bar{K}_{\varepsilon}$ in $\bar{\mathcal{H}}_{\varepsilon}$ is the
perturbation of the Kepler problem. The term $\bar{H}_{\varepsilon}$ is the
perturbation of the Hamiltonian of the PCR3BP in polar coordinates $\bar{H}_{0}$, given by
\begin{align}
\bar{H}_{0} & =\frac{1}{2}\left( \bar{R}_{2}^{2}+\left( \frac{\bar{\Phi }_{2}%
}{\bar{r}_{2}}\right) ^{2}\right) -\bar{\Phi}_{2} \notag \\
& \quad-\frac{1-\mu}{\sqrt{\bar{r}_{2}^{2}+\mu^{2}-2\mu\bar{r}_{2}\cos\left(
\phi_{2}\right) }} \notag \\
& \quad-\frac{\mu}{\sqrt{r_{2}^{2}+\left( 1-\mu\right) ^{2}+2\left(
1-\mu\right) \bar{r}_{2}\cos\left( \phi_{2}\right) }}.
\label{eqn:PCRTBP_polar}
\end{align}
\end{remark}

\begin{remark}
A customary way to deal with the  three-body problem is  to pass to  Delaunay variables  or  to Deprit variables.
However, the polar coordinates that we use are explicit, hence more convenient for our  computer assisted proof.
\end{remark}

\subsection{Geometry of the restricted three-body problem}
\label{sec:PCRTBP}
We briefly recall some of the geometric structures that organize the dynamics in the PCR3BP.
The Hamiltonian of this problem in polar coordinates is given in \eqref{eqn:PCRTBP_polar}.
%In Cartesian rotating coordinates, where the primaries
%have masses $\mu$ and $1-\mu$ ,  the distance between the primaries  is   $1$, and the period of the motion of the primaries is $2\pi$,
%the Hamiltonian is given by
%\begin{equation}\label{eq:Hamiltonian}
%\bar{H}_0=\frac{1}{2}(P_X^2+P_Y^2)+Y P_X-X P_Y-\frac{1-\mu}{r_0}-\frac{\mu}{r_1},
%\end{equation}
%where   $r_0=((X-\mu)^2+Y^2)^{1/2}$,  $r_1=((X-\mu+1)^2+Y^2)^{1/2}$,
%$X$, $Y$ are the  generalized coordinates,
%$P_X=\dot X-Y$, $P_Y=\dot Y+X$,  are the generalized momenta,
%and the symplectic form is:
%\[dP_X\wedge dX+dP_Y\wedge dY.\]
This is a $2$-degree-of-freedom Hamiltonian  and each trajectory is confined to some $3$-dimensional energy manifold $\{\bar{H}_0=h\}$ in the $4$-dimensional phase-space.

The system has 5 equilibrium points, denoted $L_1, \dotsc, L_5$.
Here we adopt the convention that $L_1$ is located between the primaries.
The equilibria $L_1,L_2,L_3$ are of saddle $\times$ center $\times$ center  linear stability type,
and the  equilibria $L_4,L_5$ are of center $\times$ center $\times$ center  linear stability type  (provided that $\mu$ is less than Routh's critical value $\mu_{\textrm{cr}}$, which is the case for our system).
A general reference for the PCR3BP  is~\cite{szebehely1967theory}.

In this paper the dynamics near the equilibrium point $L_1$ plays a crucial role.  For energy levels sufficiently close to the energy level of $L_1$, the PCR3BP possesses a family
of Lyapunov periodic orbits about $L_1$. Each Lyapunov
periodic orbit $\lambda_h$  is uniquely characterized by some fixed value of the energy $\bar{H}_0=h$.
For each periodic orbit $\lambda_h$ there exist associated stable and unstable manifolds $W^s(\lambda_h)$ and $W^u(\lambda_h)$, which are $2$-dimensional.
They are foliated by stable and unstable fibers $W^s(z)$, $W^u(z)$, $z\in\lambda_h$, respectively.
For large sets of energy values and mass parameters, it has been shown that $W^s(\lambda_h)$ and $W^u(\lambda_h)$  intersect transversally within  $\{\bar{H}_0=h\}$ (see \cite{wilczak2003heteroclinic,wilczak2005heteroclinic}).
Each such intersection represents a homoclinic point $z_h$ that  gives rise to a homoclinic orbit, which is an integral curve asymptotic in forward and backwards time to $\lambda_h$. Since a homoclinic point lies on some stable fiber $W^s(z_+)$ and also on some unstable fiber $W^u(z_-)$, for some $z_-,z_+\in\lambda_h$,   there are well defined hyperbolic stable and unstable directions at $z_h$, given by the tangent directions  at $z_h$ to $W^s(z_+)$ and  $W^u(z_-)$, respectively.

 In Section \ref{sec:coordinates} we will use these geometric objects as landmarks to construct local coordinates
 %some other objects
 for the full three-body problem (in the $6$-dimensional space).

%TCIDATA{Version=5.00.0.2606}
%TCIDATA{LaTeXparent=0,0,MMFedit.tex}

\section{Statement of the main result\label{sec:main-statement}}
Consider the full three-body problem  with initial conditions for $m_0$ and $m_1$  corresponding to the circular Kepler problem. This is  expressed by the Hamiltonian \eqref{eq:H-as-2bp+rbp-again}.
%Consider the  normalised mass of Triton in the Neptune-Triton system, which is $\mu = 0.0002089$.
Let $\Psi_t^{\varepsilon}$ stand for the flow of \eqref{eq:H-as-2bp+rbp-again}. We now state our main result and follow with comments.

\begin{theorem}\label{th:main-exact}
For the  planar full three-body problem where the mass ratio of the large bodies is $\mu = 0.0002089$ and the small body has mass $\varepsilon^2$,  there exists $\varepsilon_0$ such that,  for every $\varepsilon \in (0,\varepsilon_0]$ there exists a time $T=T(\varepsilon)$ and a point $x=x(\varepsilon)$ such that
\begin{equation}
|\bar H_{\varepsilon}(\Psi_T^{\varepsilon}(x)) - \bar H_{\varepsilon}(x)|>\Delta,
\label{eq:NT-energy-change}
\end{equation}
where $\varepsilon_0:=10^{-10}$ and $\Delta :=10^{-11}$.
\end{theorem}

Our setting corresponds to a Neptun-Triton-asteroid system, where the physical (not normalized) masses of Neptune and Triton are  $m_1=1.024 \cdot  10^{26}$ kg and $m_0 = 2.1389 \cdot 10^{22}$ kg, respectively. Theorem \ref{th:main-exact} gives   the largest physical mass  of the asteroid for which we obtain diffusion in energy  is
\[ m_2=(10^{-10})^2 \cdot 1.024 \cdot 10^{26} \, \mbox{kg} > 10^6\, \mbox{kg}.\]
%\blue Our setting corresponds to a Neptun-Triton-asteroid system, where the mass of Triton is $2.1389 \cdot 10^{22}$ kg, and  the mass of Neptune is  $1.024 \cdot  10^{26}$ kg. Theorem \ref{th:main-exact} gives that the largest physical mass  of the asteroid for which we obtain diffusion in energy is
%\[ \varepsilon_0^2 \cdot 1.024 \cdot 10^{26} \, \mbox{kg}  =(10^{-10})^2 \cdot 1.024 \cdot 10^{26} \, \mbox{kg} > 10^6\, \mbox{kg}.\]\black

The size of the energy change in (\ref{eq:NT-energy-change}) is very small, but independent of the size of the perturbation $\varepsilon$.

During the course  of the proof we will also show that the orbits we obtain, which exhibit the energy change (\ref{eq:NT-energy-change}), require physical (Earth) time $T(\varepsilon)$ between
\begin{equation}
T(\varepsilon) \in \left [\frac{\Delta}{\varepsilon}\cdot  10^3 \mbox{ years},\frac{\Delta}{\varepsilon} \cdot  25 \cdot 10^3  \mbox{ years} \right]. \label{eq:total-times}
\end{equation}
See Remark~\ref{rem:final_remark_2}.

Since $ \mathcal{\bar H}_{\varepsilon}$ is a constant of motion, from (\ref{eq:H-as-2bp+rbp-again}) we see that a change in $ \bar H_{\varepsilon}$  by $\Delta$  implies the change of $ \bar K_{\varepsilon}$ by the same amount in the opposite direction. Since $\bar H_0$ is the energy of the PCR3BP, the change in $ \bar H_{\varepsilon}$ represents a physical change of the energy of the asteroid. From (\ref{eq:K-in-bar-coordinates}) we see that a $\Delta$ change in $ \bar H_{\varepsilon}$, and thus the same order change in $\bar K_{\varepsilon}$, results in an $O(\varepsilon^2 \Delta)$ change in the Keplerian energy $K$ in the original coordinates. This means that to change the Keplerian energy $K$ by order $O(1)$ we would need the time to be of order $O(\varepsilon^{-3})$, and, moreover, the change of the energy of the restricted problem $ \bar H_{\varepsilon}$ would need to be of order $O(\varepsilon^{-2})$. We thus see that it is impossible to obtain an $O(1)$ change of $K$ for $\varepsilon$ arbitrarily close to zero, without the asteroid approaching arbitrarily close to collision. This means that our result is optimal in the sense that for an asteroid within a bounded domain away from collision,  $K$ can only change by   $O(\varepsilon^2)$. 
%TCIDATA{Version=5.00.0.2606}
%TCIDATA{LaTeXparent=0,0,MMFedit.tex}

\section{Surfaces of sections and local coordinates\label{sec:coordinates}}

\begin{figure}[tbp]
\begin{center}
\includegraphics[height=4cm]{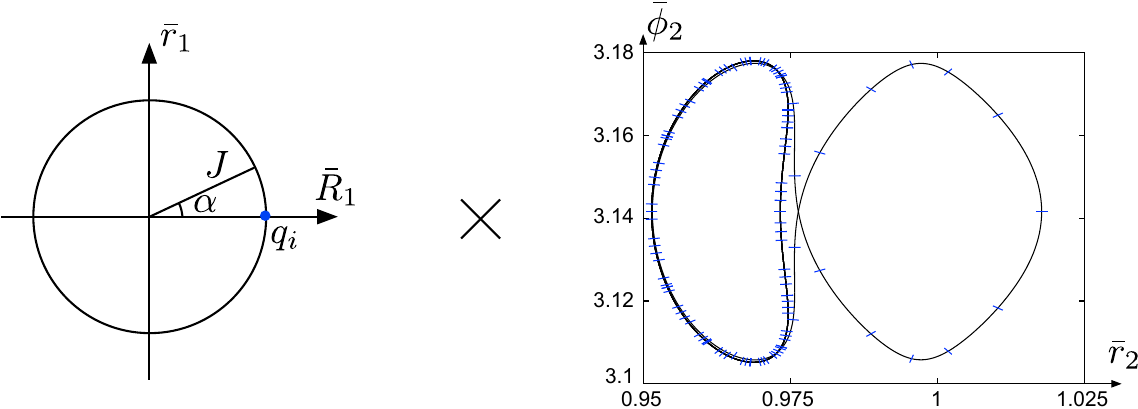}
\end{center}
\caption{On the right hand side we have a homoclinic orbit to one of the
Lyapunov orbits of the PCR3BP. The points $q_i$, for $i=0,\ldots,122$, are
positioned along the homoclinic in the $(\bar{r}_2,\bar{\phi}_2,\bar{R}_2,\bar{\Phi}_2)$ coordinates. For all the points we choose $\bar{r}_1=0$ and $
\bar{R}_1=\frac{1}{10}$.}
\label{fig:surfaces}
\end{figure}

For the proof of our result we will start with
a homoclinic orbit to a Lyapunov periodic orbit at some initial energy for the PCR3BP.
The homoclinic orbit is chosen to go around the smallest of the primaries $m_0$.
Then we will construct some surfaces of section  positioned
at a family of points $\{q_{i}\}_{i\in \{0,\ldots ,122\}}$ chosen along the homoclinic, which are
depicted in Figure \ref{fig:surfaces}. We will consider section-to-section
maps, expressed in appropriate local coordinates. The phase  space of the
system is $6$-dimensional, and by using these sections we  reduce the
dimension  to $5$.

On each section we choose coordinates so that they are aligned with the
stable, unstable and the center directions (see Section \ref{sec:PCRTBP}).
The stable and unstable directions are associated with the coordinates $(\bar{r}_{2},\bar{R}_{2},\bar{\phi}_{2},\bar{\Phi}_{2})$ of the restricted problem.
The center
directions are associated to the coordinates $(\bar{r}_{1},\bar{R}_{1})$, which we
will further express  in terms of polar coordinates in Section~\ref%
{sec:section-choices}.

On each section we also restrict to an energy level $\{\mathcal{\bar{H}}%
_{\varepsilon }=h\}$, which reduces the dimension from $5$ to $4$. We also
perform an additional local change of coordinates, which puts the
derivative of the section-to-section maps into diagonal form. This is
done in Section~\ref{sec:local-coord-on-section}.

As a result, we obtain local maps defined in $\mathbb{R}^{4}$, where the
first two coordinates are the unstable and stable directions associated with
the PCR3BP, and the second two are the radius-angle coordinates
associated with the Kepler problem.

\subsection{Choices of sections\label{sec:section-choices}}

We start by passing to polar coordinates on the $(\bar{r}_{1},\bar{R}_{1})$
plane, taking $(\alpha ,J)\in \mathbb{R}^2$ and setting
\begin{equation}
\bar{R}_{1}=J\cos (\alpha ),\qquad \bar{r}_{1}=J\sin (\alpha ).
\label{eq:to-polar-1}
\end{equation}%
(See left hand side of Figure \ref{fig:surfaces}.) This way, for $%
\varepsilon =0$ we obtain
\begin{equation}
J^{\prime }=0\qquad \alpha ^{\prime }=1.  \label{eq:dJ-dalpha}
\end{equation}%
 Note that the change of coordinates (\ref%
{eq:to-polar-1}) is not symplectic, so the formulae for $J^{\prime }$ and $%
\alpha ^{\prime }$ are derived from the vector field.

Thus, we obtain the  non-symplectic change of coordinates
\[
(\bar{r}_1, \bar{R}_1,\bar{r}_{2},\bar{\phi}_{2},\bar{R}_{2},\bar{\Phi}_{2}) 
\mapsto 
(\alpha ,J,\bar{r}_{2},\bar{\phi}_{2},\bar{R}_{2},\bar{\Phi}_{2}).
\]
%\[   (\bar{r}_1, \bar{R}_1,\bar{r}_2,\bar{R}_2, \bar{\theta}_1,\bar{\Theta}_1,\bar{\theta}_2,\bar{\Theta}_2)
%\mapsto
% (\alpha ,J,\bar{r}_2,\bar{R}_2, \bar{\theta}_1,\bar{\Theta}_1,\bar{\theta}_2,\bar{\Theta}_2).
%\]
The Hamiltonian $\mathcal{\bar{H}}_{\varepsilon }$ induces a vector field in
the coordinates $(\alpha ,J,\bar{r}_{2},\bar{\phi}_{2},\bar{R}_{2},\bar{\Phi}_{2})$, which we will denote as
$F_{\varepsilon }:\mathbb{R}^{6}\rightarrow \mathbb{R}^{6}$. Our sections are
determined by the behaviour of the system for $\varepsilon =0$.

We consider a sequence of points  $\left\{ q_{0},\ldots ,q_{122}\right\} $, where each point $q_{i}$, for $i=0,\ldots ,122$, is of the form
\begin{equation}
q_{i}=\left( 0,\frac{1}{10},\bar{r}_{2},\bar{\phi}_{2},\bar{R}_{2},\bar{\Phi}%
_{2}\right) \in \mathbb{R}^{6},  \label{eq:qi-choices}
\end{equation}%
for some $(\bar{r}_{2},\bar{R}_{2},\bar{\phi}_{2},\bar{\Phi}_{2})\in \mathbb{R}^4
$, which are along the homoclinic of the PCR3BP to a Lyapunov periodic orbit
as in Figure \ref{fig:surfaces}.
We provide our choice of the points $q_i$ in Appendix~\ref{sec:proof-K-reduced}.
 Our points are chosen
so that%
\begin{equation}
q_{0}=q_{98}=q_{122}.  \label{eq:q-start-finish}
\end{equation}%
Intuitively, we start with $q_{0}$ on a Lyapunov orbit, and choose $%
q_{1},\ldots ,q_{97}$ along the homoclinic, which leads back to $%
q_{98}=q_{0} $. The sequence of points $q_{98},$ $q_{99},\ldots ,q_{121},$ $%
q_{122}$ also starts and finishes at $q_{0}$, but the points $q_{99},\ldots
,q_{121}$ are positioned around the Lyapunov orbit. The fact that we have (%
\ref{eq:q-start-finish}) will  turn out to be important in our
construction, since it will allow us to consider the same surface of section at $q_{0}, q_{98}$ and $q_{122}$.
%We will explain the reason later on.\marginpar{Where?}

Let us use the following notation
\begin{equation*}
\pi _{2}\left( \alpha ,J,\bar{r}_{2},\bar{\phi}_{2},\bar{R}_{2},\bar{\Phi}%
_{2}\right) :=\left( 0,0,\bar{r}_{2},\bar{\phi}_{2},\bar{R}_{2},\bar{\Phi}%
_{2}\right) .
\end{equation*}
At a given point $q_{i}$ we consider a section $\Sigma _{i},$ which is
orthogonal to $\pi _{2}F_{0}(q_{i})$. In more detail,
\begin{equation*}
\Sigma _{i}:=\left\{ p:\left\langle \pi _{2}F_{0}\left( q_{i}\right)
,p-q_{i}\right\rangle =0\right\} \subset \mathbb{R}^{6},
\end{equation*}%
where $\left\langle \cdot ,\cdot \right\rangle $ is the standard scalar
product in $\mathbb{R}^{6}$. In practice, we consider a sequence of vectors $%
v_{1},\ldots ,v_{5}$ such that%
\begin{equation*}
\left\langle \pi _{2}F_{0}\left( q_{i}\right) ,v_{k}\right\rangle =0\qquad
\text{for }k=1,\ldots ,5,
\end{equation*}%
which means that
\begin{equation*}
\Sigma _{i}=q_{i}+\mathrm{span}\left\{ v_{1},\ldots ,v_{5}\right\} \subset
\mathbb{R}^{6}.
\end{equation*}
(The choice of $v_{1},\ldots ,v_{5}$ depends on the point $q_i$.)

We define a matrix $A_{i}$ by  placing the vectors $v_{1},\ldots ,v_{5}$ on the first $5$ columns of $A_{i}$, and placing  the vector $F_{0}\left( q_{i}\right) $
on the last column.  Then
\begin{equation*}
\Sigma _{i}=\left\{ q_{i}+A_{i}\left( w,0\right) :w\in \mathbb{R}^{5}\right\} .
\end{equation*}%
We can therefore identify $\Sigma _{i}$ with $\mathbb{R}^{5}$ via the
natural parameterisation of $\Sigma _{i}$ chosen as%
\begin{equation}
\mathbb{R}^{5}\ni w\mapsto q_{i}+A_{i}\left( w,0\right) \in \Sigma _{i}.
\label{eq:section-parametrisation}
\end{equation}%
We also have a natural map in the opposite direction%
\begin{equation*}
\Sigma _{i}\ni p\mapsto \pi _{w}\left( A_{i}^{-1}\left( p-q_{i}\right)
\right) \in \mathbb{R}^{5}.
\end{equation*}

Let us also introduce the notation $\mathcal{A}_{i}\in \mathbb{R}^{6\times
5}$ for the matrix whose columns coincide with the first five columns of $%
A_{i}$. Thus
\begin{equation*}
\Sigma _{i}=\left\{ q_{i}+\mathcal{A}_{i}w:w\in \mathbb{R}^{5}\right\} .
\end{equation*}

The way that we choose the vector basis $v_{1},\ldots ,v_{5}$ in practice is as follows. We
choose $v_{1}$ and $v_{2}$ to be the linear directions of the stable and
unstable fibres of the homoclinic at $q_{i}$ (see Section \ref{sec:PCRTBP}), projected onto the space
orthogonal to $\pi _{2}F_{0}\left( q_{i}\right) $. We also choose
\begin{align}
v_{3} &:=\left( 1,0,0,0,0,0\right) \in \mathbb{R}^{6},  \label{eq:alpha-v3}
\\
v_{4} &:=\left( 0,1,0,0,0,0\right) \in \mathbb{R}^{6}.  \label{eq:v4-def}
\end{align}%
Recall that our coordinates are $\left( \alpha ,J,\bar{r}_{2},\bar{R}_{2},\bar{\phi}_{2},%
\bar{\Phi}_{2}\right) $, so $v_{3}$ and $v_{4}$ are vectors
associated with $\alpha $ and $J$, respectively. They are also orthogonal to
$\pi _{2}F_{0}\left( q_{i}\right) $. We let
\begin{equation}
v_{5}:=\left( 0,0,\nabla \bar{H}_{0}\right) ,  \label{eq:v5-def}
\end{equation}%
which is also orthogonal to $\pi _{2}F_{0}\left( q_{i}\right) $.

With such choice of $v_{1},\ldots ,v_{5}$, we have the parameterisation of $%
\Sigma _{i}$ through (\ref{eq:section-parametrisation}). We have a natural
interpretation of the coordinates. For
\begin{equation*}
w=\left( w_{1},\ldots ,w_{5}\right) =\left( w_{1},w_{2},\alpha
,J,w_{5}\right) \in \mathbb{R}^{5}
\end{equation*}
the first two coordinates $w_{1}$ and $w_{2}$ are the directions of the
hyperbolic expansion and contraction, respectively. The coordinates $\left(
w_{3},w_{4}\right) =\left( \alpha ,J\right) $ are of `center'-type
coordinates, and they correspond to the coordinates  $(\bar{r}_{1},\bar{R}_{1})$ of the Kepler
problem. The last coordinate $w_{5}$ is also a `center'-type,
and is associated to $\bar{H}_{0}$.

\subsection{Choices of local coordinates on sections\label%
{sec:local-coord-on-section}}

In this subsection we introduce two additional local changes of coordinates
on our sections $\Sigma _{i}.$ The first allows us to restrict to a constant
energy level, and, at the same time ensure that one of the coordinates
corresponds to $\bar{K}_{\varepsilon }$. The second coordinate change will
allow us to have the derivatives of the section-to-section maps along the
flow to be close to diagonal.

We start with the first change of coordinates. For a fixed $\varepsilon \geq
0$ we consider an energy level%
\begin{equation*}
h:=\mathcal{\bar{H}}_{\varepsilon }\left( q_{0}\right) .
\end{equation*}%
Here we intentionally choose the energy level $h$ to be that of  the point $q_{0}$, and use the same $h$ for
all $q_{i}\in \left\{ q_{0},\ldots ,q_{122}\right\} $\footnote{%
The $q_{0},\ldots ,q_{122}$ are computed numerically, so there is no
guarantee that all the points on all the sections lie precisely on the same energy
level, due to the numerical rounding.}. We also define (see (\ref{eq:K_0-formula}), (\ref{eq:to-polar-1})
and (\ref{eq:qi-choices}))%
\begin{equation*}
\kappa _{0}:=\bar{K}_{0}\left( q_{0}\right) =\bar{K}_{0}\left( q_{i}\right) =%
\frac{1}{2}\left( \frac{1}{10}\right) ^{2}=\frac{1}{200}.
\end{equation*}

Let us focus once again on one of the $q_{i}\in \left\{ q_{0},\ldots,q_{122}\right\} $. Recall that $\mathcal{\bar{H}}_{\varepsilon }$ is a
constant of motion. Let us write $\mathrm{x}=\left( w_{1},w_{2},\alpha,I
\right) \in \mathbb{R}^{4}$ and implicitly define a function $\omega _{i}:%
\mathbb{R}^{4}\rightarrow \mathbb{R}^{5}$ such that for $\mathrm{x}=\left(
w_{1},w_{2},\alpha ,I\right) \in \mathbb{R}^{4}$ the function $\omega_i$ satisfies
\begin{eqnarray}
\mathcal{\bar{H}}_{\varepsilon }\left( q_{i}+\mathcal{A}_{i}\left( \omega
_{i}\left( \mathrm{x}\right) \right) \right)  &=&h,
\label{eq:omega-implicit} \\
\bar{K}_{\varepsilon }\left( q_{i}+\mathcal{A}_{i}\left( \omega _{i}\left(
\mathrm{x}\right) \right) \right)  &=&\kappa _{0}+I,
\label{eq:I-is-Keps-kappa}
\end{eqnarray}

The choice of local coordinates $\mathrm{x}=\left(
w_{1},w_{2},\alpha ,I\right)$ on $\Sigma_i$ with
\begin{equation}
\left( w_{1},w_{2},\alpha ,J,w_{5}\right) =\omega _{i}\left(
w_{1},w_{2},\alpha ,I\right)   \label{eq:from-I-to-J}
\end{equation}%
ensures that the points $q_{i}+\mathcal{A}_{i}\left( \omega _{i}\left(
\mathrm{x}\right) \right) \in \Sigma _{i}$ will all have the energy level
equal to $\mathcal{\bar{H}}_{\varepsilon }=h$. This is ensured by (\ref%
{eq:omega-implicit}). Moreover, from (\ref{eq:I-is-Keps-kappa}) we see that
the change of coordinates (\ref{eq:from-I-to-J}) ensures that $I=\bar{K}%
_{\varepsilon }-\kappa _{0}$.

For a fixed $\mathrm{x}=\left( w_{1},w_{2},\alpha ,I\right) \in \mathbb{R}%
^{4}$ the value of $\omega_{i}\left( \mathrm{x}\right) $ can be found by
solving numerically
\begin{equation}
\mathcal{G}_{i}\left( \mathrm{x},y;\varepsilon\right) =0  \label{eq:G-newton-function}
\end{equation}%
 for $y\in \mathbb{R}^{2}$ for the function $\mathcal{G}_{i}:\mathbb{R}^{4}\times \mathbb{R}^{2}\rightarrow \mathbb{R}^{2}$
defined as
\begin{eqnarray*}
\mathcal{G}_{i}\left( \mathrm{x},y; {\varepsilon }\right)  &=&\mathcal{G}_i%
\left( \left( w_{1},w_{2},\alpha ,I\right) ,\left(
y_{4},y_{5}\right); {\varepsilon } \right)  \\
&:=&\left(
\begin{array}{c}
\mathcal{\bar{H}}_{\varepsilon }\left( q_{i}+\mathcal{A}_{i}\left(
w_{1},w_{2},\alpha ,y_{4},y_{5}\right) \right) -h \\
\bar{K}_{\varepsilon }\left( q_{i}+\mathcal{A}_{i}\left( w_{1},w_{2},\alpha
,y_{4},y_{5}\right) \right) -\kappa _{0}-I%
\end{array}%
\right) .
\end{eqnarray*}%
Then for $y=\left( y_{4},y_{5}\right) $ which solves $\mathcal{G}_{i}
\left( \mathrm{x},y;\varepsilon\right) =0$ we have
\begin{equation*}
\omega_{i} \left( \mathrm{x}\right) =\omega_{i} \left( w_{1},w_{2},\alpha ,I\right)
=\left( w_{1},w_{2},\alpha ,y_{4},y_{5}\right) .
\end{equation*}%
(Here we write $y_{4},y_{5}$ to emphasize the connection between these variables
and  the vectors $v_{4},v_{5}$ from (\ref{eq:v4-def}--\ref{eq:v5-def}).)

The equation
(\ref{eq:G-newton-function}) can be solved by  using an  interval arithmetic-based  Newton method, which can give rigorous  interval bounds for $\omega_{i} \left( \mathrm{x}\right) $.

Moreover, from (\ref{eq:omega-implicit}) we
obtain the derivative of $\omega_{i}  \left( \mathrm{x}\right) $ as%
\begin{equation*}
D\omega_{i}  =\left(
\begin{array}{c}
\begin{array}{cccc}
1 & 0 & 0 & 0 \\
0 & 1 & 0 & 0 \\
0 & 0 & 1 & 0%
\end{array}
\\
-\left( \frac{\partial \mathcal{G}_{i} }{\partial y}\right) ^{-1}\frac{\partial
\mathcal{G}_{i} }{\partial \mathrm{x}}%
\end{array}%
\right) \in \mathbb{R}^{5\times 4}.
\end{equation*}

We now describe our second change of coordinates, which will allow us to
obtain that the derivatives of the section-to-section maps expressed in the
local coordinates will be close to diagonal. Consider a matrix $B_{i}\in
\mathbb{R}^{4\times 4}$ and a vector $\mathbf{w}_{i}\in \mathbb{R}^{4}$ and
define new coordinates $\mathbf{x}=\left( u,s,\alpha,I \right)$ on $\Sigma_i$ by
\[(w_1,w_2,\alpha,I)= B_i (u,s,\alpha,I) + \varepsilon \mathbf{w}_i=B_i \mathbf{x} + \varepsilon \mathbf{w}_i.\]
We always choose matrices $B_{i}$ and vectors $\mathbf{w}_{i}$ of the form%
\begin{eqnarray}
B_{i} &=&\left(
\begin{array}{cccc}
1 & 0 & b_{13}^{i} & b_{14}^{i} \\
0 & 1 & b_{23}^{i} & b_{24}^{i} \\
0 & 0 & 1 & 0 \\
0 & 0 & 0 & 1%
\end{array}%
\right) ,  \label{eq:B-form} \\
\mathbf{w}_{i} &=&\left( \mathbf{w}_{u}^{i},\mathbf{w}_{s}^{i},0,0\right) ,
\label{eq:w-form}
\end{eqnarray}%
which ensures that
\begin{equation}
\pi _{\alpha ,I}\left( B_{i}\mathbf{x}+\varepsilon \mathbf{w}_{i}\right)
=\pi _{\alpha ,I}\mathbf{x}.  \label{eq:B-I-preserved}
\end{equation}

Thus, we obtain a parametrization
\[\gamma_i^\varepsilon :\mathbb{R}^4 \to\Sigma_i\] of $\Sigma_i$ in terms of the new coordinates $\mathbf{x}$, which is of the form
\begin{equation}
\gamma_i^\varepsilon (\mathbf{x})=\gamma_i^\varepsilon (u,s,\alpha,I)=(\alpha,J,\bar r_2, \bar \phi_2,\bar R_2,\bar \Phi_2),\label{eq:total-coord-change}
\end{equation}
where the formula for $\gamma_i^\varepsilon$ is
\[\gamma_i^\varepsilon(\mathbf{x}) =q_i + \mathcal{A}_i \,\omega_i( B_i \mathbf{x} + \varepsilon \mathbf{w}_i).\]

\black Our setup has an important property stated in the lemma and the remarks below.

\begin{lemma}
\label{lem:I-on-section} For a point on a section $\Sigma _{i}$
expressed in the local coordinates
$\mathbf{x}=\left( u,s,\alpha,I \right)$, we have%
\begin{equation*}
\bar{K}_{\varepsilon }\left( \gamma _{i}^{\varepsilon }\left( \mathbf{x}%
\right) \right) =\kappa _{0}+I.
\end{equation*}
\end{lemma}

\begin{proof}
From (\ref{eq:omega-implicit}) and (\ref{eq:B-I-preserved}) we have
\begin{equation*}
\bar{K}_{\varepsilon }\left( \gamma _{i}^{\varepsilon }\left( \mathbf{x}%
\right) \right) =\kappa _{0}+\pi _{I}\left( B_{i}\mathbf{x}+\varepsilon
\mathbf{w}_{i}\right) =\kappa _{0}+I,
\end{equation*}%
as required.
\end{proof}

\begin{remark}
Lemma \ref{lem:I-on-section} ensures that if we construct an orbit which
has an increase $\Delta $ in the local coordinate $I$, then such orbit
increases by $\Delta $ in $\bar{K}_{\varepsilon }$.
\end{remark}

\begin{remark}
All our coordinate changes preserve the coordinate $\alpha $, meaning that
the  $\alpha $ in the local coordinates $(u,s,\alpha,I)$ is the same as the  $\alpha $ introduced in (\ref{eq:to-polar-1}).
\end{remark}

We now discuss the particular choices we make of the coefficients of the
matrices $B_{i}$ and the vectors $\mathbf{w}_{i}$ which enter into our
changes of coordinates (\ref{eq:total-coord-change}). The points $q_{0},q_{1},\ldots ,q_{122}$  that we consider
lie along a homoclinic orbit of the PCR3BP.

The  section-to-section maps along the flow  $\Psi^\varepsilon_{t}$  are defined by
\begin{equation}\label{eqn:section-to-section}
\begin{split}
\mathcal{P}_{i}^{\varepsilon }:&\Sigma _{i-1}\rightarrow \Sigma _{i},\\
\mathcal{P}_{i}^{\varepsilon }(x)&=\Psi^\varepsilon_{\tau_i^\varepsilon(x)}(x),
\end{split}
\end{equation}
where
\begin{equation}\label{eqn:tau_i}
  \tau_i^\varepsilon(x) =\textrm{ the shortest time $t$  to take $x\in \Sigma_{i-1}$ to $\Psi^\varepsilon_{t}  (x)\in\Sigma _{i}$.}
\end{equation}

Our choice of local coordinates implies that for the section-to-section maps $\mathcal{P}_{i}^{\varepsilon }$
we have
\begin{equation}
\pi _{2}\mathcal{P}_{i}^{\varepsilon }\left( q_{i-1}\right) \approx \pi
_{2}q_{i}\qquad \text{for }i=1,\ldots ,122.  \label{eq:poinc-qi}
\end{equation}%
At each point $q_{i}$ we have a coordinate change $\gamma _{i}^{\varepsilon }
$ of the form (\ref{eq:total-coord-change}). We can then define local maps
\begin{equation}
\begin{split}
f_{i}^{\varepsilon }:&\mathbb{R}^{4}\rightarrow \mathbb{R}^{4},\\
f_{i}^{\varepsilon }:&=\left( \gamma _{i}^{\varepsilon }\right) ^{-1}\circ
\mathcal{P}_{i}^{\varepsilon }\circ \gamma _{i-1}^{\varepsilon }.
\end{split}
 \label{eq:local-maps-i}
\end{equation}%
We choose the coefficients of $B_{i}$ and the vectors $\mathbf{w}_{i}$ for
the coordinate changes so that we obtain
\begin{equation*}
\begin{array}{c}
\frac{\partial \left( \pi _{u,s}f_{i}^{\varepsilon =0}\right) }{\partial
\left( I,\alpha \right) }\approx 0 \\
\frac{\partial \left( \pi _{u,s}f_{i}^{\varepsilon =0}\right) }{\partial
\varepsilon }\approx 0%
\end{array}%
\qquad \text{for }i=1,\ldots ,122.
\end{equation*}%
With such choice the derivatives of the local maps take the form%
\begin{equation}
Df_{i}^{\varepsilon }\approx \left(
\begin{array}{llll}
\lambda _{i} & 0 & 0 & 0 \\
0 & 1/\lambda _{i} & 0 & 0 \\
0 & 0 & 1 & \kappa _{i} \\
0 & 0 & 0 & 1%
\end{array}%
\right) ,  \label{eq:good-alignment-of-local-maps}
\end{equation}%
where $\lambda _{i}$ and $1/\lambda _{i}$ are the expansion and contraction
rates, respectively, and $\kappa _{i}$ is a twist coefficient. From the choices (\ref%
{eq:total-coord-change}) (\ref{eq:B-form}), (\ref{eq:w-form}) (\ref%
{eq:B-I-preserved}) we have
\begin{equation*}
\gamma _{i}^{\varepsilon =0}\left( 0\right) =q_{i}+\mathcal{A}_{i}\left(
\omega \left( 0\right) \right) =q_{i}.
\end{equation*}

To summarize, through our change of coordinates we are able to achieve the
following goals:

\begin{enumerate}
\item Our system is reduced from the $6$-dimensional flow
to $4$-dimensional maps.

\item We ensure that the local coordinates are mapped to appropriate
surfaces of sections, restricted to a constant energy level $\mathcal{\bar H}%
_{\varepsilon}=h$.

\item A change in the local coordinate $I$ corresponds precisely to a
change in $\bar K_{\varepsilon}$.

\item The local maps \eqref{eq:local-maps-i} are well aligned with the expanding/contracting/center
coordinates, so that we obtain (\ref{eq:good-alignment-of-local-maps}).
\end{enumerate}

Finally we ensure that%
\begin{equation}
\Sigma _{0}=\Sigma _{98}=\Sigma _{122}\qquad \text{and\qquad }\gamma
_{0}^{\varepsilon }=\gamma _{98}^{\varepsilon }=\gamma _{122}^{\varepsilon },
\label{eq:sections-same}
\end{equation}%
by choosing%
\begin{eqnarray}
q_{0} &=&q_{98}=q_{122},\qquad A_{0}=A_{98}=A_{122},  \label{eq:q-A-B-same-1}
\\
B_{0} &=&B_{98}=B_{122},\qquad w_{0}=w_{98}=w_{122}.  \label{eq:q-A-B-same-2}
\end{eqnarray}%
This will play an important role in the proof of Theorem \ref{th:main-intro} in section \ref{sec:proof},  since it ensures that we have the same local coordinates on the section $\Sigma _{0}=\Sigma _{98}=\Sigma _{122}$.
%later on.
%\black \marginpar{Where?}

%TCIDATA{Version=5.00.0.2606}
%TCIDATA{LaTeXparent=0,0,MMFedit.tex}

\section{Tools for establishing diffusion}
\label{sec:tools}
In the previous section we have defined local maps (\ref{eq:local-maps-i})
which are expressed in local coordinates%
\begin{equation}
\mathbf{x}=\left( u,s,\alpha ,I\right) \in \mathbb{R}^{4},
\label{eq:sec-coord-2}
\end{equation}%
introduced in (\ref{eq:total-coord-change}). We have ensured in Lemma \ref%
{lem:I-on-section} that $I=\bar{K}_{\varepsilon }-\kappa _{0}$, where $%
\kappa _{0}$ is a constant.

Since for $\varepsilon =0$ the energy of the Kepler problem $\bar{K}_{0}$ is
a constant of motion, by Lemma \ref{lem:I-on-section} we see that%
\begin{equation*}
\pi _{I}f_{i}^{\varepsilon =0}\left( \mathbf{x}\right) =\pi _{I}%
\mathbf{x}, \qquad \mbox{for }i=1,\ldots, 122.
\end{equation*}

 In this section we review the diffusion mechanism from \cite{MR4544807}, which  establishes the existence of
orbits  along which the action coordinate $I$ changes by an amount  independent
of the size of the perturbation $\varepsilon $.
In Section~\ref{sec:proof} we verify this diffusion mechanism  for our system.
By Lemma~\ref{lem:I-on-section} the change in $I$ will imply that we have the same change in $\bar{K}_{\varepsilon }$.

The method is based on appropriate topological alignment of sets by the
local maps $f_{i}^{\varepsilon }$, which we refer to as correct
alignment of windows, and on appropriate cone conditions. These tools are introduced in Section~\ref{sec:cover-cc}.
In Section \ref{sec:diffusion-mechanism} we recall Theorem~\ref{th:diffusion-mechanism}, which is our main tool for establishing
diffusion in $I$.

\subsection{Correctly aligned windows and cone conditions\label{sec:cover-cc}}

In this section we recall two notions, correctly aligned windows and cone
conditions, which are needed to introduce the main tool for establishing
Arnold diffusion. For more details, see \cite{GideaZ04a,zgliczy2009covering,MR4544807}.

We start with correctly aligned windows.

We shall write $B^{n}\left( z,r\right) $ to denote a ball of radius $r$ in $%
\mathbb{R}^{n},$ centered at $z\in \mathbb{R}^{n}$. We shall write $B^{n}$
for a unit ball centered at zero in $\mathbb{R}^{n}$. For a set $A\subset
\mathbb{R}^{n}$ we shall denote its closure as $\bar{A}$ and its boundary as
$\partial A$.

A window in $\mathbb{R}^{n}=\mathbb{R}^{n_{u}}\times \mathbb{R}^{n_{cs}}$ is a set of the form%
\begin{equation*}
N=\bar{B}^{n_{u}}\left( z_{1},r_{1}\right) \times \bar{B}^{n_{cs}%
}\left( z_{2},r_{2}\right) ,
\end{equation*}%
where $z_{1}\in \mathbb{R}^{n_{u}}$, $z_{2}\in \mathbb{R}^{n_{cs}}$
and $r_{1},r_{2}\in \mathbb{R}$. The notation $n_{u}$ here stands for
`topologically unstable' dimension and $n_{cs}$ stands for `topologically stable' dimension.
%Unlike in the hyperbolic case, the topologically unstable directions do  not need to be uniformly expanding, and the topologically stable directions do  not need to be uniformly contracting.
%\marginpar{I commented out a sentence here} % MC: The reason is that possibly people might not be familiar with correctly aligned windows, in which case the the hyperbolic case and this comment will not tell them much and might create confusion.

\begin{figure}
\begin{center}
\includegraphics[height=2cm]{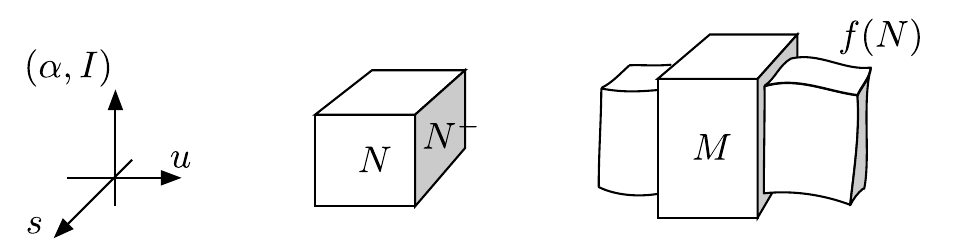}
\end{center}
\caption{ An example of correctly aligned windows $N\overset{f}{\implies }M$. The exit sets $N^-$ and $M^-$ are depicted in grey. The entry sets $N^+$ and $M^+$ consist of the white sides of the cubes $N$ and $M$, respectively. Here $cs=(s,\alpha,I)$, where $s$ is contracting and $\alpha,I$ are center coordinates. \label{fig:covering}}
\end{figure}

We introduce an
`exit set' and an `entry set' of $N$, defined respectively as%
\begin{equation*}
N^{-}=\partial \bar{B}^{n_{u}}\left( z_{1},r_{1}\right) \times \bar{B}^{n_{cs}}\left( z_{2},r_{2}\right) ,\qquad N^{+}=\bar{B}^{n_{u}%
}\left( z_{1},r_{1}\right) \times \partial \bar{B}^{n_{cs}}\left(
z_{2},r_{2}\right) .
\end{equation*}%
For a point $z\in \mathbb{R}^{n_{u}}\times \mathbb{R}^{n_{cs}}$ we
shall write $z=\left( u,cs\right)$, to distinguish its coordinates. See Figure \ref{fig:covering}.

\begin{definition}
Let $N$ and $M$ be two windows in $\mathbb{R}^{n}=\mathbb{R}^{n_{u}%
}\times \mathbb{R}^{n_{cs}}$. Let $f:N\rightarrow \mathbb{R}^{n}$ be a
continuous mapping. We say that $N$ is correctly aligned with $M$, and write%
\begin{equation*}
N\overset{f}{\implies }M
\end{equation*}

if the following conditions are satisfied:

\begin{enumerate}
\item There exists a homotopy $\chi :\left[ 0,1\right] \times N\rightarrow
\mathbb{R}^{n_{u}}\times \mathbb{R}^{n_{cs}}$ such that%
\begin{equation*}
\chi \left( 0,\cdot \right) =f\left( \cdot \right) ,\qquad \chi \left( \left[
0,1\right] ,N^{-}\right) \cap M=\emptyset ,\qquad \chi \left( \left[ 0,1%
\right] ,N\right) \cap M^+=\emptyset ,
\end{equation*}

\item There exists a linear map $A:\mathbb{R}^{n_{u}}\rightarrow
\mathbb{R}^{n_{u}}$ such that%
\begin{eqnarray*}
\chi \left( 1,u,cs\right) &=&\left( Au,0\right) \qquad \text{for all }\left(
u,cs\right) \in N\subset \mathbb{R}^{n_{u}}\times \mathbb{R}^{\mathrm{cs%
}}, \\
A\left( \partial B^{n_{u}}\right) &\subset &\mathbb{R}^{n_{u}%
}\setminus \bar{B}^{n_{u}}.
\end{eqnarray*}
\end{enumerate}
\end{definition}

\begin{remark}
The above definition of correctly aligned windows, as well as the cone conditions
introduced below,  work in arbitrary dimension. In the setting of this
paper we have $n_{u}=1$ and $n_{cs}=3.$ This is associated with
our choice of local coordinates $\left( u,s,\alpha ,I\right) \in \mathbb{R}^{%
n_{u}}\times \mathbb{R}^{n_{cs}}=\mathbb{R}^{4}$, and we will consider  $u$  as the topologically unstable
unstable coordinate, $cs=(s,\alpha ,I)$ as the `topologically stable' coordinate.
The latter notation is justified  since $cs$  includes both the stable coordinate $s$ and the center coordinates $(\alpha ,I)$.
\end{remark}

We now introduce  cone conditions. Consider a function $Q:\mathbb{R}%
^{n_{u}}\times \mathbb{R}^{n_{cs}}\rightarrow \mathbb{R}$%
\begin{equation}
Q\left( u,cs\right) =\left\Vert u\right\Vert _{n_{u}}^{2}-\left\Vert
cs\right\Vert _{n_{cs}}^{2},  \label{eq:cone-Q}
\end{equation}%
where $\left\Vert \cdot \right\Vert _{n_{u}}$ and $\left\Vert \cdot
\right\Vert _{n_{cs}}$ stand for some norms in $\mathbb{R}^{n_{u}}$
and $\mathbb{R}^{n_{cs}}$, respectively. For a point $z\in \mathbb{R}^{%
n_{u}}\times \mathbb{R}^{n_{cs}}$ we define the $Q$-cone at $z$ as
the set $\{z^{\prime }:Q\left( z-z^{\prime }\right) >0\}$. We shall refer to
the function $Q$ as a cone.

\begin{figure}
\begin{center}
\includegraphics[height=2cm]{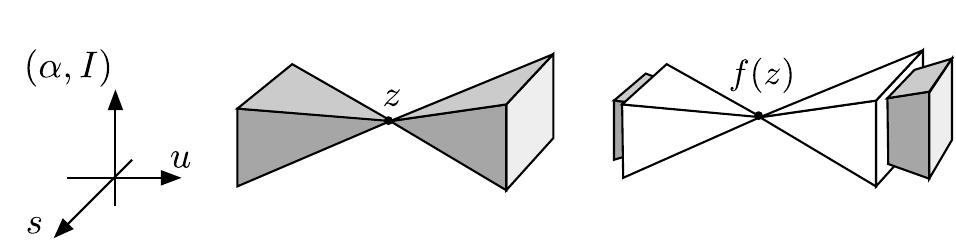}
\end{center}
\caption{Cone condition. A gray cone at $z$ is mapped inside a white cone at $f(z)$.  \label{fig:cone}}
\end{figure}

\begin{definition}
Let $Q_{1},Q_{2}:\mathbb{R}^{n_{u}}\times \mathbb{R}^{n_{cs}%
}\rightarrow \mathbb{R}$ be cones of the form (\ref{eq:cone-Q}). We say that
a function $f:N\rightarrow \mathbb{R}^{n_{u}}\times \mathbb{R}^{\mathrm{%
cs}}$ satisfies cone conditions from $Q_{1}$ to $Q_{2}$ if for $z,z^{\prime
}\in N$ (see Figure \ref{fig:cone})
\begin{equation*}
Q_{1}\left( z-z^{\prime }\right) >0\qquad \text{implies\qquad }Q_{2}\left(
f\left( z\right) -f\left( z^{\prime }\right) \right) >0.
\end{equation*}
\end{definition}

\begin{definition} \label{def:windows_cones}
Given two windows, $N_{1}$ and $N_{2}$,   two cones $Q_{1}$ and $Q_{2}$ on $N_{1}$ and $N_{2}$, respectively,
and  a continuous map $f:N_{1}\rightarrow \mathbb{R}^{n_{u}}\times
\mathbb{R}^{n_{cs}}$,  we say that $N_{1}$ is correctly aligned with $N_{2}$ with cone conditions,
if  $N_1$ is correctly aligned with $N_2$ under $f$,
and moreover $f$ satisfies cone conditions from $Q_{1}$ to $Q_{2}$.   We
denote this by
\begin{equation*}
\left( N_{1},Q_{1}\right) \overset{f}{\implies }\left( N_{2},Q_{2}\right) .
\end{equation*}
\end{definition}

\begin{remark}\label{rem:horizontal_disc}
A related notion is that of  a horizontal disc, which is a graph of a continuous map $h:\bar{B}^{n_{u}}\to\mathbb{R}^{n_{u}}\times \mathbb{R}^{n_{cs}}$.
A horizontal disc satisfies a $Q$-cone condition if $x\neq x'$ in $\bar{B}^{n_{u}}$ implies $Q(h(x)-h(x'))>0$.
If $h$ is a horizontal disc in the window $N_{1}$, and $N_{1}$ is correctly aligned with $N_{2}$
under  $f$ satisfying cone conditions from $Q_{1}$ to $Q_{2}$,  then for any horizontal disc $h$  in $N_1$ satisfying a
 a $Q_1$-cone condition, the image of $h$ under $f$ contains a horizontal disc   in $N_2$ satisfying
 a $Q_2$-cone condition.
For details see \cite{MR4544807}.
\end{remark}

\subsection{Diffusion mechanism\label{sec:diffusion-mechanism}}

Recall that the local coordinate $I$ is equal to $\bar{K}_{\varepsilon
}-\kappa _{0},$ where $\bar{K}_{\varepsilon }$ is the Keplerian part of the
Hamiltonian of the three-body problem (\ref{eq:H-as-2bp+rbp}) and $\kappa
_{0}$ is a constant. For $\varepsilon =0$ the $I$ is an integral of motion.
In Section \ref{sec:proof} will show that for each $\varepsilon >0$ sufficiently small there are orbits along which the action $I$ changes by an amount independent of $\varepsilon $. Since the total energy $\mathcal{%
\bar{H}}_{\varepsilon }$ is preserved, from (\ref{eq:H-as-2bp+rbp}) we see
that the change  in $I$ will imply a transfer of energy from the Keplerian
part $\bar{K}_{\varepsilon }$ to the restricted three-body problem part $%
\bar{H}_{\varepsilon }$.

In Section \ref{sec:section-choices} we have introduced a sequence of
surfaces of sections $\Sigma _{0},\ldots ,\Sigma _{122}.$ We have ensured
that%
\begin{equation*}
\Sigma _{0}=\Sigma _{98}=\Sigma _{122}.
\end{equation*}%
We have equipped $\Sigma _{i}\cap \{ \mathcal{\bar{H}}_{\varepsilon
}=h\} $, for $i=0,\ldots ,122,$ with local coordinates $(u,s,\alpha ,I)
$ introduced by (\ref{eq:total-coord-change}). We have ensured in (\ref%
{eq:sections-same}) that the local coordinates on the sections $\Sigma
_{0},\Sigma _{98}$ and $\Sigma _{122}$ are the same.

Since $\alpha $ is an angle we will identify angles which are equal mod $%
2\pi $.

Let us consider a set $\mathbf{S}\subset \mathbb{R}^{4}$ of the form%
\begin{equation}
\mathbf{S}=\bar{B}^{n_{u}}\times \bar{B}^{n_{s}}\times S_{\alpha
}\times \mathbb{R}  \label{eq:strip-S-def}
\end{equation}%
where $S_{\alpha }$ is an interval%
\begin{equation*}
S_{\alpha }=\left[ \alpha _{1},\alpha _{2}\right] .
\end{equation*}%
We will refer to the set $\mathbf{S}$ as the `strip'.
In our case we have $n_{u}=n_{s}=1$.

We will now consider a finite collection of windows $\left\{ N_{\ell ,0}\right\}
_{\ell \in L}$ on the section $\Sigma _{0}$, where $L$ is a finite set of
indexes. Each window $N_{\ell ,0}\subset \mathbb{R}^{4}$, for $\ell \in L$,
is expressed in the local coordinates $(u,s,\alpha ,I)$ on $\Sigma _{0}\cap
\{\mathcal{\bar{H}}_{\varepsilon }=h\}$.

In Section~\ref{sec:proof} we will choose  $\# L=9\cdot 10^4$.
We will also provide the reason  for choosing  a family of $L$ initial windows, which is required by the implementation of validated numerics.

Each $N_{\ell ,0},$ for $\ell \in L$, will be an initial window of a
sequence of correctly aligned windows. We will assume that each window $%
N_{\ell ,k}$ involved in such sequence is expressed in the local coordinates
$\left( u,s,\alpha ,I\right) $ on $\Sigma _{k}\cap \{\mathcal{\bar{H}}%
_{\varepsilon }=h\}.$ Each sequence will end with a window $N_{\ell ,k_{\ell
}}$, where $k_{\ell }$ can depend on the choice of $\ell \in L$. We will
assume that the initial window $N_{\ell ,0}$ and the final window $N_{\ell
,k_{\ell }}$ are on the section $\Sigma _{0}=\Sigma _{98}=\Sigma _{122}$,
meaning that%
\begin{equation*}
k_{\ell }\in \left\{ 98,\,122\right\} .
\end{equation*}%
We will also equip the sections $\Sigma _{0}$, $\Sigma _{98}$ and $\Sigma
_{122}$ with the same cone $Q$.

We now make the above more precise by giving a definition of a connecting sequence \cite{MR4544807}.

\begin{definition}
\label{def:connection-sequence} (Connecting sequence)

\begin{enumerate}
\item A connecting sequence consists of a sequence of windows%
\begin{equation*}
N_{\ell ,0},N_{\ell ,1},\ldots ,N_{\ell ,k_{\ell }},
\end{equation*}%
and a sequence of cones%
\begin{equation*}
Q_{\ell ,0},Q_{\ell ,1},\ldots ,Q_{\ell ,k_{\ell }},
\end{equation*}%
such that we have the following correct alignment of windows and cone
conditions%
\begin{equation*}
\left( N_{\ell ,0},Q_{\ell ,0}\right) \overset{f_{1}^{\varepsilon }}{%
\implies }\left( N_{\ell ,1},Q_{\ell ,1}\right) \overset{f_{2}^{\varepsilon }%
}{\implies }\ldots \overset{f_{k_{\ell }}^{\varepsilon }}{\implies }\left(
N_{\ell ,k_{\ell }},Q_{\ell ,k_{\ell }}\right) .
\end{equation*}

\item The initial cone $Q_{\ell ,0}$ and the final cone $Q_{\ell ,k_{\ell }}$
are equal,%
\begin{equation}
Q_{\ell ,0}=Q_{\ell ,k_{\ell }}=Q.  \label{eq:Q-introduced}
\end{equation}%
(The cone $Q$ is independent of $\ell \in L.$)

\item For the initial window $N_{\ell ,0}$ and the final window $N_{\ell
,k_{\ell }}$ we have%
\begin{equation*}
N_{\ell ,0},N_{\ell ,k_{\ell }}\subset \mathbf{S\qquad }\text{and\qquad }\pi
_{u,s}N_{\ell ,0}=\pi _{u,s}N_{\ell ,k_{\ell }}=\pi _{u,s}\mathbf{S}=\bar{B}%
^{n_{u}}\times \bar{B}^{n_{s}}.
\end{equation*}
\end{enumerate}
\end{definition}

\begin{definition}
We shall say that a point $x\in N_{\ell ,0}$ passes through a connecting
sequence whenever%
\begin{equation*}
f_{i}^{\varepsilon }\circ \ldots \circ f_{1}^{\varepsilon }\left( x\right)
\in N_{\ell ,i},\qquad \text{for all }i=1,\ldots ,k_{\ell }.
\end{equation*}%
The time needed for such a point  $x$ to pass through the connecting sequence is
\begin{equation}
\tau _{1}^{\varepsilon }\left( x\right) +\tau _{2}^{\varepsilon }\left(
f_{1}(x)\right) +\ldots +\tau _{k_{\ell }}^{\varepsilon }\left( f_{k_{\ell
}-1}\circ \ldots \circ f_{1}(x)\right)
\label{eqn:time_needed}
\end{equation}%
where $\tau_i$ is the time needed to go from one section to another \eqref{eqn:tau_i}.
\end{definition}

Recall that the sections restricted to the constant energy level are locally
homeomorphic to $\mathbb{R}^{4}$ and that on these sections we have local
coordinates $\left( u,s,\alpha ,I\right) $. Let us assume that the cone $Q$
from (\ref{eq:Q-introduced}) is of the form%
\begin{equation}
Q\left( u,s,\alpha ,I\right) =u^{2}-\left( \max \left\{ \frac{1}{a_{s}}%
\left\vert s\right\vert ,\frac{1}{a_{\alpha }}\left\vert \alpha \right\vert ,%
\frac{1}{a_{I}}\left\vert I\right\vert \right\} \right) ^{2}
\label{eq:cone-Q-form}
\end{equation}%
where $a_{s},a_{I},a_{\alpha }>0$ are constants.

The following theorem will be our main tool for establishing Arnold
diffusion. The hypothesis of this theorem is that we have appropriate connecting sequences relative to the section-to-section maps $f^\varepsilon_i$,  which satisfy certain conditions. The conclusion of the theorem is that there are diffusing orbits for the flow $\Psi _{t}^{\varepsilon}$ of   \eqref{eq:H-as-2bp+rbp-again}.

\begin{theorem}
\cite[Theorem 4.2]{MR4544807}\label{th:diffusion-mechanism} Assume that:

\begin{enumerate}
\item[i.] For each $\ell \in L$ there we have a connecting sequence
satisfying
\begin{equation*}
\left( N_{\ell ,0},Q_{\ell ,0}\right) \overset{f_{1}^{\varepsilon }}{%
\implies }\left( N_{\ell ,1},Q_{\ell ,1}\right) \overset{f_{2}^{\varepsilon }%
}{\implies }\ldots \overset{f_{k_{\ell }}^{\varepsilon }}{\implies }\left(
N_{\ell ,k_{\ell }},Q_{\ell ,k_{\ell }}\right) .
\end{equation*}

\item[ii.] The projections of the initial windows $N_{\ell ,0}$ onto the $%
\left( I,\alpha \right) $-coordinates covers $\left[ 0,1\right] \times
S_{\alpha }$, meaning that%
\begin{equation*}
\bigcup_{\ell \in L}\pi _{\alpha,I }N_{\ell ,0}= S_{\alpha }\times \left[ 0,1%
\right] .
\end{equation*}%
(Recall that the interval $S_{\alpha }$ was used to define the strip $%
\mathbf{S}$ in (\ref{eq:strip-S-def}).)

\item[iii.] Whenever $N_{\ell ,0}\cap N_{\ell ^{\prime },0}\neq \emptyset $,
for every $I^{\ast },\alpha ^{\ast }\in N_{\ell ,0}\cap N_{\ell ^{\prime
},0} $ the multidimensional rectangle%
\begin{equation*}
\bar{B}^{n_{u}}\times \bar{B}^{n_{s}}\times \left( \left[ I^{\ast
}-a_{I},I^{\ast }+a_{I}\right] \cap \left[ 0,1\right] \right) \times \left( %
\left[ \alpha ^{\ast }-a_{\alpha },\alpha ^{\ast }+a_{\alpha }\right] \cap
S_{\alpha }\right)
\end{equation*}%
is contained in $N_{\ell ,0}$ or $N_{\ell ^{\prime },0}.$

\item[iv.] For every $x$ which passes through the connecting sequence we have%
\begin{equation*}
\varepsilon C>\pi _{I}\left( f_{k_{\ell }}^{\varepsilon }\circ \ldots \circ
f_{1}^{\varepsilon }\left( x\right) -x\right) >c\varepsilon ,
\end{equation*}%
for some $C>c>0$ which are independent of $\ell .$
\end{enumerate}

Then for every $\varepsilon >0$ there exists a point $x=x\left( \varepsilon
\right) \in \mathbf{S}$ and a time $T=T\left( \varepsilon \right) $ such that%
\begin{equation*}
\pi _{I}x=0\qquad \text{and\qquad }\pi _{I}\Psi _{T}^{\varepsilon }\left(
x\right) >1.
\end{equation*}

Moreover, the diffusion time $T\left( \varepsilon \right) $ is bounded by
\begin{equation}
T=T\left( \varepsilon \right) \in [(C\varepsilon)^{-1} \delta,
(c\varepsilon)^{-1}\rho],   \label{eq:diffusion-time-from-mechanism}
\end{equation}
where $\rho >0$  and  $\delta>0$ are upper and lower bounds, respectively, for the time needed to pass through each connecting sequence (see \eqref{eqn:time_needed}).
\end{theorem}

\begin{remark}\label{rem:connecting}
In Theorem \ref{th:diffusion-mechanism} condition $i$  says that we have a finite family of correctly aligned windows with cone conditions so that the initial and the final window for each connecting sequence lands in the same strip. Moreover, by condition $ii$  the projections of the initial windows onto the center coordinates  cover the center component of the strip. The correct alignment of windows with cone conditions implies the existence of horizontal discs (see Remark \ref{rem:horizontal_disc} ) that propagate from the first  to the last window in each connecting sequence. Condition $iii$ implies that a disc that has been propagated
along a connecting sequence can be propagated again along some other connecting sequence.
Condition $iv$ says that points that get propagated along each connecting sequence get their action coordinate $I$ increases by $O(\varepsilon)$.
If these conditions are fulfilled, Theorem \ref{th:diffusion-mechanism} implies the existence of diffusing orbits along which the action coordinate $I$ increases by $O(1)$. For details, see \cite{MR4544807}.
\end{remark}

\begin{remark} An important feature of Theorem \ref{th:diffusion-mechanism} is that it requires the validation of a finite number of connecting sequences, each of a finite length. In all, this requires a finite number of estimates, which can be validated by a computer using interval arithmetic. The length of the trajectories which exhibits the macroscopic changes in $I$ is of order $O(1/\varepsilon)$, which becomes arbitrarily long as $\varepsilon$ goes to zero. The mechanism outlined in Remark \ref{rem:connecting} allows us to shadow such arbitrarily long trajectories even though we validate a finite number of conditions. 
\end{remark}

%TCIDATA{Version=5.00.0.2606}
%TCIDATA{LaTeXparent=0,0,MMFedit.tex}

\section{Proof of the main result}
\label{sec:proof}
In this section we provide the proof of our main Theorem \ref{th:main-intro},
which has been formulated in a more detailed form in Theorem \ref%
{th:main-exact}. We start with Section \ref{sec:numerics}, where we provide
numerical evidence and intuition that our setup will lead to the verification
of the assumptions of Theorem \ref{th:diffusion-mechanism}. Then in Section %
\ref{sec:the-proof} we will apply Theorem \ref{th:diffusion-mechanism} to
give a computer assisted proof of Theorem \ref{th:main-intro}.

For our computer assisted proof we use the CAPD library \cite{MR4283203}. It provides a rigorous interval arithmetic integrator which relies on the Lohner's method \cite{MR1930946}. The library can evaluate rigorous enclosures of section-to-section maps  along the flow, together with bounds on their derivatives \cite{MR4395996}.

\subsection{The numerical properties of the homoclinic orbit\label%
{sec:numerics}}

The local coordinates $(u,s,\alpha ,I)$ introduced in Section \ref%
{sec:coordinates} have led us to local maps which satisfy (\ref%
{eq:good-alignment-of-local-maps}). This means that we will have contraction
and expansion on the the coordinates $u$ and $s$, respectively, which is
needed for the correct alignment of windows. The coordinates $\alpha $ and $%
I $ are neither contracting or expanding. We will treat them as
topologically stable by enlarging the successive windows in the $\alpha $
and $I$ components.

The first issue  in this section is to construct
connecting sequences along the homoclinic so that they
start and finish in the same strip $\mathbf{S}$. The main difficulty is how to  position the initial and final windows of a connecting sequence
along the coordinate $\alpha $.
The second issue is construct
connecting sequences  such that, under the perturbation $\varepsilon >0$, the trajectories
passing through the connecting sequences will increase in $I$. In this
section we provide numerical insights to these two issues.

\begin{figure}[tbp]
\begin{center}
\includegraphics[height=4cm]{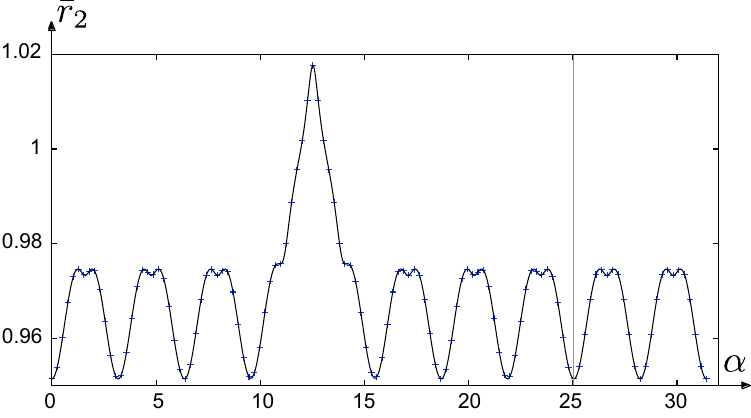}
\end{center}
\caption{The homoclinic orbit from Figure \protect\ref{fig:surfaces} plotted
in $\protect\alpha ,\bar{r}_{2}$ coordinates.}
\label{fig:alpha}
\end{figure}

The first aspect we address is the change in the coordinate $\alpha $ along
the homoclinic. In terms of the coordinate $\alpha $, an important feature
for us is that we have (\ref{eq:q-start-finish}), i.e.
\begin{equation*}
q_{0}=q_{98}=q_{122}
\end{equation*}%
and that along $\alpha $ the iterates of the section-to-section maps $%
\mathcal{P}_{i}^{\varepsilon }:\Sigma _{i-1}\rightarrow \Sigma _{i}$ of the
point $q_{0}$
\begin{equation}
\alpha _{i}:=\pi _{\alpha }\mathcal{P}_{i}^{\varepsilon =0}\circ \ldots
\circ \mathcal{P}_{1}^{\varepsilon =0}\left( q_{0}\right)
\label{eq:alpha_i-def}
\end{equation}%
behave as follows%
\begin{equation}
\begin{array}{lll}
\alpha _{98}=25.05192798, & \qquad \qquad & \alpha _{98}
=-0.080813244 \quad \text{ mod }2\pi, \\
\alpha _{122}=31.4174449, & \qquad \qquad & \alpha _{122}
=0.001518367 \quad \text{ mod }2\pi.%
\end{array}
\label{eq:alpha-shifts}
\end{equation}
This property is depicted in Figure \ref{fig:alpha} where we plot the
homoclinic orbit in coordinates $\alpha ,\bar{r}_{2}$. (Compare with Figure %
\ref{fig:surfaces}.) The points marked with blue crosses indicate where the trajectory crosses
the consecutive sections $\Sigma _{1},\ldots ,\Sigma _{122}$. The gray
vertical line indicates where the trajectory crosses the section $\Sigma
_{98}.$

The important observation  is that for the point $\mathcal{P}_{98}^{\varepsilon
=0}\circ \ldots \circ \mathcal{P}_{1}^{\varepsilon =0}\left( q_{0}\right) $
the angle $\alpha $ (mod $2\pi $) is decreased with respect to $\pi _{\alpha
}q_{0}=0$, and that for the point $\mathcal{P}_{122}^{\varepsilon =0}\circ
\ldots \circ \mathcal{P}_{1}^{\varepsilon =0}\left( q_{0}\right) $ the  angle $\alpha $  (mod $2\pi $)
is increased with respect $\pi _{\alpha }q_{0}=0$.
This will ensure that for $\varepsilon $ close to zero we will have
\begin{eqnarray*}
\pi _{\alpha }f_{98}^{\varepsilon }\circ \ldots \circ f_{1}^{\varepsilon
}\left( x\right)  &<&\pi _{\alpha }x, \\
\pi _{\alpha }f_{122}^{\varepsilon }\circ \ldots \circ f_{1}^{\varepsilon
}\left( x\right)  &>&\pi _{\alpha }x.
\end{eqnarray*}

This will be useful for the validation of the assumptions of Theorem \ref%
{th:diffusion-mechanism}. We will choose a strip $\mathbf{S}$ of the from (%
\ref{eq:strip-S-def}) with%
\begin{equation*}
S_{\alpha }=\left[ 0,0.0825\right] .
\end{equation*}%
Such strip is wide enough so that from (\ref{eq:alpha-shifts}) we see that
for $x\in \mathbf{S}$ with $\pi _{\alpha }x$ from zero up to $%
0.0825-0.00152=0.08098$ the $\pi _{\alpha }f_{122}^{\varepsilon }\circ
\ldots \circ f_{1}^{\varepsilon }\left( x\right) $ will land back in $%
S_{\alpha }$. Similarly, for $\pi _{\alpha }x$ from $0.08082$ to $0.0825$
the $\pi _{\alpha }f_{98}^{\varepsilon }\circ \ldots \circ
f_{1}^{\varepsilon }\left( x\right) $ will also land back in $S_{\alpha }$.
Since the two intervals%
\begin{equation}
S_{\alpha }^{l}:=\left[ 0,0.08098\right] \qquad \text{and\qquad }S_{\alpha
}^{r}:=\left[ 0.08082,0.0825\right]   \label{eq:Sl-Sr-def}
\end{equation}%
overlap and their union is $S_{\alpha }$ we see that this  ensures that
for every point $x\in \mathbf{S}$ we are able to choose either the
sequence $f_{98}^{\varepsilon },\ldots ,f_{1}^{\varepsilon }$ or the
sequence $f_{122}^{\varepsilon },\ldots ,f_{1}^{\varepsilon }$ of local maps
to return back to $S_{\alpha }$ along the coordinate $\alpha $. This will be
essential  to verify property \textit{3.} from the Definition \ref%
{def:connection-sequence} for our connecting sequences constructed for the
validation of assumptions of Theorem \ref{th:diffusion-mechanism}. Namely,
the behaviour of the local maps along $\alpha $ will allow us to choose the
windows so that $N_{\ell ,k_{\ell }}\subset \mathbf{S}$.

We now turn to the issue of verifying that  the coordinate $I$  increases along connecting sequences.
The
points $q_{0},\ldots ,q_{122}$ have the property that a trajectory starting
from points close to $q_{0}$ will increase in $I$ under a perturbation $%
\varepsilon >0$. To describe this in more detail, first let us observe that our choice of the local coordinates (\ref%
{eq:total-coord-change}) implies that%
\begin{equation*}
\gamma _{i}^{\varepsilon =0}\left( 0\right) =q_{i}.
\end{equation*}%
From (\ref{eq:alpha_i-def}) we also have%
\begin{equation*}
f_{i}^{\varepsilon =0}\circ \ldots \circ f_{0}^{\varepsilon =0}\left(
0\right) =\left( 0,\alpha _{i}\right) ,
\end{equation*}%
where  $(0,\alpha _{i})\in \mathbb{R}^{3}\times \mathbb{R}$.

We verify that
that the points $q_{0},\ldots ,q_{122}$ have the property that for%
\begin{equation*}
c_{i}:=\frac{d}{d\varepsilon }\left( \pi _{I}f_{i}^{\varepsilon }\left(
0,\alpha _{i}\right) \right) |_{\varepsilon =0}\qquad \text{for }i=1,\ldots
,122,
\end{equation*}%
we have
\begin{equation}
\min \left( c_{1}+\ldots +c_{98},\,c_{1}+\ldots +c_{122}\right) >c>0
\label{eq:c-bound-for-qi}
\end{equation}%
for some constant $c\in \mathbb{R}$. The property (\ref{eq:c-bound-for-qi})  follows from our careful choice of the homoclinic
orbit.

The reason why (\ref{eq:c-bound-for-qi}) is important for us is that since $%
\pi _{I}f_{i}^{\varepsilon =0}\left( x\right) =\pi _{I}x$, for small $%
\varepsilon $ and for $x$ close to $\left( 0,\alpha _{i}\right) $ we have%
\begin{equation*}
\pi _{I}f_{i}^{\varepsilon }\left( x\right) -\pi _{I}x\approx \varepsilon
c_{i}.
\end{equation*}%
This means that our points $q_{0},\ldots ,q_{122}$ have the property that an orbit starting close to $q_{0}$ will increase in $I$ under the
composition of the local maps $f_{k_{\ell }}^{\varepsilon }\circ \ldots
\circ f_{1}^{\varepsilon }$ for $k_{\ell }\in \{98,122\}$. This will be
important for the verification of  the assumption \textit{iv} from Theorem \ref%
{th:diffusion-mechanism}.

\subsection{The proof of Theorem \protect\ref{th:main-intro}\label%
{sec:the-proof}}

Our proof is divided into steps, which are reflected in Algorithm \ref{algorithm} which is presented at the end of the proof. 

 \subsubsection*{Step 0: Choice of the local maps and the strip}\black We choose the sequences $q_{i}\in \mathbb{R}^{6},$ $A_{i}\in \mathbb{R}%
^{6\times 6},$ $B_{i}\in \mathbb{R}^{4\times 4}$ and $w_{i}\in \mathbb{R}^{4}
$, which define the sections $\Sigma _{i}$ and local coordinates on sections
$\Sigma _{i}$, for $i=0,\ldots ,122.$ We ensure by (\ref{eq:q-A-B-same-1}--%
\ref{eq:q-A-B-same-2}) that the sections%
\begin{equation*}
\Sigma _{0}=\Sigma _{98}=\Sigma _{122}
\end{equation*}%
are equipped with the same local coordinates. The particular choices of $%
q_{i},$ $A_{i},$ $B_{i},$ $w_{i}$ are given in the data tables in \cite{code}.
We emphasize here that regardless of the choice of the particular values
for these points and matrices, the computer code rigorously ensures that for a point $%
\mathbf{x}\in \mathbb{R}^{4}$ in any local coordinates, we have%
\begin{equation*}
\bar{K}_{\varepsilon }\left( \gamma _{i}^{\varepsilon }\left( \mathbf{x}%
\right) \right) =\kappa _{0}+\pi _{I}\left( \mathbf{x}\right) .
\end{equation*}%
(See Lemma \ref{lem:I-on-section}.) If our computer assisted proof
establishes the changes in the local coordinate $I$, then these changes have
to correspond to the actual energy transfer between $\bar{K}_{\varepsilon }$
and $\bar{H}_{\varepsilon }$ in the full three body problem.

We choose
\begin{equation}
r:=10^{-9} \label{eq:r-choice}
\end{equation}%
and the following strip on $\Sigma _{0}$
\begin{equation*}
\mathbf{S:}=\bar{B}^{n_{u}}\left( r\right) \times \bar{B}^{n_{s}%
}\left( r\right) \times S_{\alpha }\times \mathbb{R}
\end{equation*}%
with%
\begin{equation*}
\bar{B}^{n_{u}}\left( r\right) =\bar{B}^{n_{s}}\left( r\right)
:=[-r,r]\qquad \text{and\qquad }S_{\alpha }=\left[ 0,0.0825\right] .
\end{equation*} 

\subsubsection*{Step 1: Choice of initial windows and initial cones} 
We subdivide the interval $S_{\alpha }=\left[ 0,0.0825\right] $ into $ 9\cdot 10^{4}$ overlapping sub intervals $\mathbf{a}_{\ell }$, for $\ell \in L:=\{1,\ldots,  9\cdot
10^{4}\}$, with an overlap
which ensures the conditions\textit{\footnote{%
In our proof we consider a small interval $\left[ 0,10^{-11}\right] $ of $I$
instead of $\left[ 0,1\right] $, so we require that $\bigcup_{\ell }\pi
_{\alpha ,I}N_{\ell ,0}=S_{\alpha }\cap \left[ 0,10^{-11}\right] $. Also,
since we consider $\bar{B}^{\mathrm{u}}\left( r\right) =\bar{B}^{\mathrm{s}%
}\left( r\right) =[-r,r]$ and not $\left[ -1,1\right] ,$ for \textit{iii}
from Theorem \ref{th:diffusion-mechanism} we need that%
\begin{equation*}
\bar{B}^{\mathrm{u}}\times \bar{B}^{\mathrm{s}}\times \left( \left[ I^{\ast
}-ra_{I},I^{\ast }+ra_{I}\right] \cap \left[ 0,10^{-11}\right] \right)
\times \left( \left[ \alpha ^{\ast }-ra_{\alpha },\alpha ^{\ast }+ra_{\alpha
}\right] \cap S_{\alpha }\right)
\end{equation*}%
is contained in $N_{\ell ,0}$ or $N_{\ell ^{\prime },0}.$ See \cite[Remark
4.4]{MR4544807}.}} \textit{ii }and \textit{iii} from Theorem \ref%
{th:diffusion-mechanism}. 

For an interval $\mathbf{a}_{\ell }\subset
S_{\alpha }$ we define
\begin{equation}
N_{\ell ,0}:=\bar{B}^{\mathrm{u}}\left( r\right) \times \bar{B}^{\mathrm{s}%
}\left( r\right) \times \mathbf{a}_{\ell }\times \left[ 0,10^{-11}\right] \label{eq:N0-choice}
\end{equation}%
and choose (see (\ref{eq:Sl-Sr-def}))
\begin{equation}
k_{\ell }=\left\{
\begin{array}{ll}
122 & \qquad \text{if }\mathbf{a}_{\ell }\subset S_{\alpha }^{l}, \\
98 & \qquad \text{if }\mathbf{a}_{\ell }\subset S_{\alpha }^{r}.%
\end{array}%
\right.   \label{eq:choice-of-length}
\end{equation}

Our particular choice of the parameters defining
the cone $Q$ on $\Sigma _{0}$ are $a_{s}=a_{\alpha }=a_{I}=10^{-2}$ (see (%
\ref{eq:cone-Q-form})), giving
\begin{eqnarray}
Q\left( u,s,\alpha ,I\right)  &=& u^2 - \left( \max \left\{ \frac{1%
}{a_{s}}\left\vert s\right\vert ,\frac{1}{a_{\alpha }}\left\vert \alpha
\right\vert ,\frac{1}{a_{I}}\left\vert I\right\vert \right\}\right)^2 \notag \\
&=& u^2 -\left(10^{2}\,\max \left\{ \left\vert s\right\vert
,\left\vert \alpha \right\vert ,\left\vert I\right\vert \right\} \right)^2 . \label{eq:Q-choice}
\end{eqnarray}% 

\subsubsection*{Step 2: Validation of a connecting sequence} 
Starting with each initial window $N_{\ell ,0}$ we construct a connecting
sequence as follows.

Each set $N_{\ell ,i}$ for $i=1,\ldots ,k_{\ell }$ is constructed to be
\begin{equation*}
N_{\ell ,i}=\bar{B}^{n_{u}}\left( r\right) \times \bar{B}^{n_{s}%
}\left( r_{\ell ,i}\right) \times \mathbf{a}_{\ell ,i}\times \mathbf{I}_{\ell ,i}
\end{equation*}%
where the $r_{\ell ,i}$, $\mathbf{a}_{\ell ,i}$, $\mathbf{I}_{\ell ,i}$ are closed
intervals in $\mathbb{R}$ chosen to satisfy
\begin{equation}
\pi _{s,\alpha ,I}\left( f_{i}\left( N_{\ell ,i-1}\right) \right) \subset \bar{B}^{n_{s}}\left( r_{\ell ,i}\right) \times
\mathbf{a}_{\ell ,i}\times \mathbf{I}_{\ell ,i}.  \label{eq:a-I-choice}
\end{equation}
The choice of $r_{\ell ,i}$, $\mathbf{a}_{\ell ,i},\mathbf{I}_{\ell
,i}$ is done automatically by our code. For $i=k_{\ell}$ we ensure though that $r_{\ell ,k_{\ell}}=r$ with $r$ from (\ref{eq:r-choice}), which means that 
\begin{equation}
\pi_s N_{\ell,k_{\ell}} = \pi_s \mathbf{S} . \label{eg:last-window-in-S-on-s}
\end{equation} 
Once we compute the interval enclosure of $f_{i}\left( N_{\ell ,i-1}\right) $ the intervals are chosen to
be large enough to ensure (\ref{eq:a-I-choice}). 

Note that on the $u$ the sets $\pi _{u}N_{\ell ,i}=\bar{B}^{n_{u}}\left( r\right)$ are always chosen to
be of the same size, with $r$ from (\ref{eq:r-choice}). In particular, this means that 
\begin{equation}
\pi_u N_{\ell,k_{\ell}} = \pi_u \mathbf{S} . \label{eg:last-window-in-S-on-u}
\end{equation} 
 The computer assisted proof validates that for every $\varepsilon \in \left[
0,10^{-10}\right] $ \begin{equation*}
N_{\ell  ,0}\overset{f_{1}^{\varepsilon }}{\implies }%
N_{\ell ,1}\overset{f_{2}^{\varepsilon }}{\implies }\ldots \overset{%
f_{k_{\ell }}^{\varepsilon }}{\implies }N_{\ell ,k_{\ell }}.
\end{equation*}
We thus validate  a sequence of coverings of lengths $98$ or $%
122$, where the length depends on the choice of $\mathbf{a}_{\ell }$; see (%
\ref{eq:choice-of-length}).

At the same time our computer program ensures that $f_{i}^{\varepsilon }$, for $i=1,\ldots,k_{\ell}$,
satisfy cone conditions from $Q_{\ell,i-1}$ to $Q_{\ell,i}$, with the choice of the initial cone on $\Sigma_0$
\[Q_{\ell,0}=Q,\] 
with $Q$ as in (\ref{eq:Q-choice}). (The cones $Q_{\ell,i}$ on intermediate sections $\Sigma_i$, for $i=1,\ldots, k_{\ell}$ are chosen automatically by our program.)

\subsubsection*{Step 3: Checking that we return to the strip } 
For each $\ell $ we validate
that%
\begin{equation}
\pi _{\alpha }N_{\ell ,k_{\ell }}\subset S_{\alpha } \label{eq:final-alpha}
\end{equation}%
which together with (\ref{eg:last-window-in-S-on-s}) and (\ref{eg:last-window-in-S-on-u}) implies that%
\begin{equation*}
N_{\ell ,k_{\ell }}\subset \mathbf{S.}
\end{equation*}%
This is possible due to the behaviour of the local maps along $\alpha $,
described in section \ref{sec:numerics}.

\subsubsection*{Step 4: Checking the final cone condition } 

Our program automatically computes the sequence of cones $Q_{\ell,i}$ so that the local maps $f_i^{\varepsilon}$ satisfy cone conditions from $Q_{\ell,i-1}$ to $Q_{\ell,i}$, for $i=1,\ldots,k_{\ell}$. We ensure that the final cone $Q_{\ell,k_{\ell}}$ is contained in $Q$, which ensures that $f_{k_{\ell}}^{\varepsilon}$ satisfies cone conditions from $Q_{\ell,k_{\ell}-1}$ to $Q$.

\subsubsection*{Step 5: Computing the change in $I$ along a connecting sequence } 

Our code validates that for each $\ell$ and every $i\in \{1,\ldots,k_{\ell}\}$, $\varepsilon\in[  0,10^{-10}]$, $q\in\gamma_{i-1}%
^{\varepsilon}\left(  N_{\ell,i-1}\right)  $ we have%
\begin{equation}
c_{\ell,i}<\frac{d}{d\varepsilon}\left(  K_{\varepsilon}\left(  \mathcal{P}%
_{i}^{\varepsilon}\left(  q\right)  \right)  -K_{\varepsilon}\left(  q\right)
\right)  <C_{\ell
,i}.\label{eq:c_i-bounds}%
\end{equation}
Since for $\varepsilon=0$ the $K_{0}$ is a constant of motion,
\[ K_{\varepsilon=0}\left(  \mathcal{P}_{i}^{\varepsilon=0}\left(
q\right)  \right)  -K_{\varepsilon=0}\left(  q\right)   =0
\]
so we see that
\[
K_{\varepsilon}\left(  \mathcal{P}_{i}^{\varepsilon}\left(  q\right)  \right)
-K_{\varepsilon}\left(  q\right)  =\varepsilon\int_{0}^{1}\frac{\partial
}{\partial\varepsilon}\left(  K_{s\varepsilon}\left(  \mathcal{P}%
_{i}^{s \varepsilon}\left(  q\right)  \right)  -K_{s\varepsilon}\left(
q\right)  \right)  ds.
\]
Hence, from (\ref{eq:c_i-bounds}) it follows that for  $q\in\gamma
_{i-1}^{\varepsilon}\left(  N_{\ell,i-1}\right)  $
\begin{equation}
\varepsilon c_{\ell,i}<K_{\varepsilon}\left(  \mathcal{P}_{i}^{\varepsilon
}\left(  q\right)  \right)  -K_{\varepsilon}\left(  q\right)  <\varepsilon
C_{\ell,i}.\label{eq:ci-bounds-2}%
\end{equation}
Above inequality holds for $q=\gamma_{i-1}^{\varepsilon}\left(
x\right)  $, for all $x\in N_{\ell,i-1}$. Since $I=K_{\varepsilon}%
-\kappa_{0}$, we see that for $x\in N_{\ell,i-1}$ and $q=\gamma_{i-1}%
^{\varepsilon}\left(  x\right)  $
\[
\pi_{I}f_{i}^{\varepsilon}\left(  x\right)  -\pi_{I}x=K_{\varepsilon}\left(
\mathcal{P}_{i}^{\varepsilon}\left(  \gamma_{i-1}^{\varepsilon}\left(
x\right)  \right)  \right)  -K_{\varepsilon}\left(  \gamma_{i-1}^{\varepsilon
}\left(  x\right)  \right)  =K_{\varepsilon}\left(  \mathcal{P}_{i}%
^{\varepsilon}\left(  q\right)  \right)  -K_{\varepsilon}\left(  q\right)  .
\]
This means that (\ref{eq:ci-bounds-2}) ensures that for every $x\in
N_{\ell,i-1}$%
\[
\pi_{I}f_{i}^{\varepsilon}\left(  x\right)  -\pi_{I}x\in\left[  \varepsilon
c_{\ell,i},\varepsilon C_{\ell,i}\right]  ,
\]
hence%
\[
\varepsilon\sum_{i=1}^{k_{\ell}}C_{\ell,i}\geq\pi_{I}f_{k_{\ell}}%
^{\varepsilon}\circ\ldots\circ f_{1}^{\varepsilon}\left(  x\right)  -\pi
_{I}x\geq\varepsilon\sum_{i=1}^{k_{\ell}}c_{\ell,i}.
\]
The $c_{\ell,i}$ and $C_{\ell,i}$, for $i=1,\ldots,k_{\ell}$, are computed
automatically by our computer program from
(\ref{eq:c_i-bounds}). \black

%Our code validates that for every $\ell $ we have

%\begin{equation}
%c_{\ell ,i}<\frac{d}{d\varepsilon }\pi _{I}f_{i}^{\varepsilon }\left(N_{\ell ,i-1}\right) _{\mid\varepsilon \in \left[ 0,10^{-10}\right] }<C_{\ell,i}. \label{eq:c_i-bounds}
%\end{equation}%
%Since $\pi _{I}f_{i}^{\varepsilon =0}\left( x\right) =\pi _{I}x$, this ensures that for every $x\in N_{\ell ,i-1}$%
%\begin{equation*}
%\pi _{I}f_{i}^{\varepsilon }\left( x\right) -\pi _{I}x=\int_{0}^{1}\frac{d}{ds}\pi _{I}f_{i}^{s\varepsilon }\left( x\right) ds\in \left[ \varepsilon c_{\ell ,i},\varepsilon C_{\ell ,i}\right] ,
%\end{equation*}%
%hence%
%\begin{equation*}
%\varepsilon \sum_{i=1}^{k_{\ell }}C_{\ell ,i}\geq \pi _{I}f_{k_{\ell}}^{\varepsilon }\circ \ldots \circ f_{1}^{\varepsilon }\left( x\right) -\pi_{I}x\geq \varepsilon \sum_{i=1}^{k_{\ell }}c_{\ell ,i}.
%\end{equation*}
%The $c_{\ell ,i}$ and $C_{\ell ,i}$, for $i=1,\ldots,k_{\ell}$, are computed automatically by our computer program from (\ref{eq:c_i-bounds}).

 \subsubsection*{Step 6: Validating the global bound on the change in $I$ } 
Finally, we validate that for every $\ell \in L$ we
obtain%
\begin{equation}
C>\sum_{i=1}^{k_{\ell }}C_{\ell ,i}>\sum_{i=1}^{k_{\ell }}c_{\ell ,i}>c, \label{eq:final-c}
\end{equation}%
with%
\begin{equation}
%C=3.995\cdot 10^{-4}\qquad \text{and\qquad }c=1.9\cdot 10^{-5}.
 C= 3.9888 \cdot 10^{-4}\qquad \text{and\qquad }c=2.2695\cdot 10^{-5}. \label{eq:final-bounds}
\end{equation}

 \subsubsection*{Summary of steps in the computer-assisted proof} 
The steps of the proof can be summarised as follows. For every $\ell \in L$
\begin{enumerate}
\item choose $N_{\ell,0}$ and $Q_{\ell,0}=Q$ according to (\ref{eq:N0-choice}) and (\ref{eq:Q-choice}), respectively;
\item validate
        			 \[(N_{\ell,i-1},Q_{\ell,i-1})\overset{f_i^{\varepsilon}}{\implies} (N_{\ell,i},Q_{\ell,i})\qquad \mbox{for } i=1,\ldots k_{\ell},\]
for all $\varepsilon \in [0,10^{-10}]$;
\item validate (\ref{eq:final-alpha});
\item validate that the cones given by $Q_{\ell,k_{\ell}}$ are contained in the cones given by $Q$;
\item validate (\ref{eq:c_i-bounds}); (The $c_{\ell ,i}$ and $C_{\ell ,i}$ are computed by our program.)
\item validate (\ref{eq:final-c}).
\end{enumerate}
%The resulting algorithm for the computer assisted proof is presented below. 
We can see that the process is trivially parallelizable since the computations for different $\ell$ can be performed independently.

This way we  validate  that assumptions \textit{i--iv} of Theorem \ref%
{th:diffusion-mechanism} are satisfied, which implies that the claim of
Theorem \ref{th:main-intro} is true.

%\begin{algorithm}
%\caption{The algorithm for our computer assisted proof \label{algorithm}}
%\begin{algorithmic}[1]
%\Statex \textbf{for} $\ell \in L$ \textbf{do}
%	\State \quad Choose $N_{\ell,0}$ and $Q_{\ell,0}=Q$ according to (\ref{eq:N0-choice}) and (\ref{eq:Q-choice}), respectively.
%	\Statex  \quad \textbf{for} $i=1$ to $k_{\ell}$ \textbf{do}
%		\State \quad \quad Validate
%        			 \[(N_{\ell,i-1},Q_{\ell,i-1})\overset{f_i^{\varepsilon}}{\implies} (N_{\ell,i},Q_{\ell,i})\]
%		\quad \quad for all $\varepsilon \in [0,10^{-10}]$. (The $N_{\ell,i}$ and $Q_{\ell,i}$ are computed by our program.)
%\Statex  \quad \textbf{end for}
%	\State \quad  Validate (\ref{eq:final-alpha}).
%	\State \quad  Validate that the cones given by $Q_{\ell,k_{\ell}}$ are contained in the cones
%	\Statex \quad  given by $Q$.
%	\State \quad  Validate (\ref{eq:c_i-bounds}). (The $c_{\ell ,i}$ and $C_{\ell ,i}$ are computed by our program.)
%\Statex \textbf{end for}
%\State  Validate (\ref{eq:final-c}).
%\end{algorithmic}
%\end{algorithm}

\begin{remark} \label{rem:final_remark_2} The physical times from (\ref{eq:total-times}) follow from the fact that it
takes $5.877$ Earth days for Triton to revolve around Neptune and the fact
that the longest integration time from $\Sigma _{0}$ to $\Sigma _{122}$ in
our construction was under $\rho =31.5$ revolutions of Triton around
Neptune, and the shortest integration time from $\Sigma _{0}$ to $\Sigma
_{98}$ was over $\delta =25$ revolutions of Triton around Neptune; see (\ref{eq:alpha-shifts}). The bound
(\ref{eq:diffusion-time-from-mechanism}) from Theorem \ref%
{th:diffusion-mechanism} is for an energy change equal to $1$, but here we
change by $\Delta $, hence the scaling by $\Delta $ in (\ref{eq:total-times}%
). Combining (\ref{eq:final-bounds}) with (\ref{eq:diffusion-time-from-mechanism}) we
obtain (\ref{eq:total-times}).
\end{remark}

\begin{remark}\label{rem:final_remark_1} The most difficult part of the computer assisted proof is to obtain $c>0$.
To obtain sharp enough estimates for  (\ref{eq:c_i-bounds}), we
needed to consider very small sets. This forced us in particular to take
small intervals $\varepsilon \in \left[ 0,10^{-10}\right] $ and $I\in \left[
0,10^{-11}\right] .$ We also have to place many intermediate sections since
to use small sets the points $q_i$ on the homoclinic orbit need to depart and return to a
small neighbourhood of the Lyapunov orbit.

To work with small sets we also had to subdivide the strip $\mathbf{S}$
along $S_{\alpha }$ to $9\cdot10^{4}$ small initial sets $N_{\ell ,0}$. For
each initial set we had to validate a connecting sequence and bonds for (\ref{eq:final-c}). Due to this the
computation time for our computer assisted proof is substantial, equal to 26 
days on a single thread. We have performed the computations on the Athena
supercomputer at Cyfronet AGH. We have used a tiny fraction of Athena's
power. We have used its single node, running 48 parallel tasks on it. (Our
computation on Athena took under $13$  hours). 

Even though the computational time is substantial, the computer assisted proof could be performed on an average configuration desktop (say with 8 cores and 16 threads) over a weekend.
\end{remark}

\begin{remark} \label{rem:final_remark_3} We have quite a bit of flexibility in the choice of the vectors and matrices
defining the sections and local coordinates. Also, the numbers $98$ and $122$ of
the intermediate sections are quite arbitrary. Our particular choices
follow from a careful numerical exploration of the system, so that we
diagonalise the dynamics well in the local coordinates, leading to (\ref%
{eq:good-alignment-of-local-maps}). For other choices it is quite likely
that the computer assisted part of the proof would also establish the
result. In particular, if one was to marginally modify our choices, then for
sufficiently small modification the computer assisted proof would also be
successful. If one would choose these vectors and matrices in an
irrational/incorrect way, the computer program would simply report that it
cannot establish the result.
\end{remark}

\begin{figure}
	\begin{center}
		\includegraphics[height=2cm]{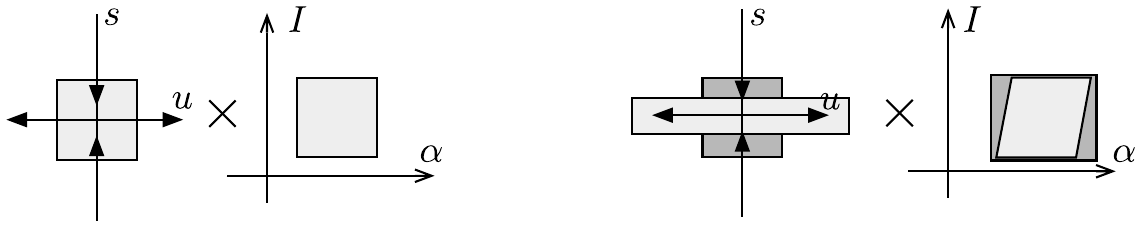}
	\end{center}
	\caption{Overcoming the wrapping effect. Our sets are cartesian products of cubes in the $(u,s)$ and in the $(\alpha,I)$ coordinates. The initial set is depicted in light grey. After integrating we can choose the second set, in darker grey, to be of the same size in the $(u,s)$ coordinates. On coordinates $(\alpha,I)$ we need the second set to be slightly larger, to enclose the image of the initial set.\label{fig:wrapping-control}}
\end{figure}

\begin{remark}

A number of numerical challenges were involved in our computer assisted proof:

{\bf Finding the domain in which we have the transfer of energies.} Our system is Hamiltonian. If we have a variable $I$, which is an integral of motion prior to perturbation, then after a Hamiltonian perturbation, due to the exact symplecticity of the maps, the domain is split into two regions: one where $I$ grows, and the other where $I$ decreases. Trajectories of the perturbed system travel through such regions and typically $I$ remains roughly at the same level. We needed to have a construction in which we can ensure that the trajectories of the perturbed system remain within a neighbourhood in which $I$, on average, is increasing. Finding such neighbourhood required a careful numerical investigation of the system.

%{\bf Explicit formulae for the vector field.} We need explicit formulae for the vector field to be able to use an interval arithmetic integrator for the ODE. (We cannot for instance use the Delaunay variables or the Deprit variables in our approach.) We have adopted the method from \cite{MR1694376} to make all the reductions explicit.

\begin{figure}
	\begin{center}
		\includegraphics[height=3cm]{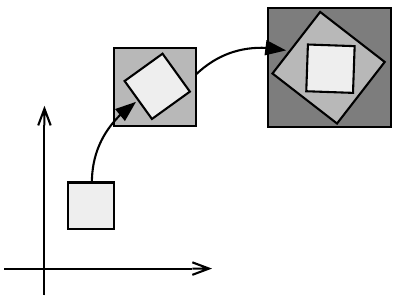}
	\end{center}
	\caption{After an interval arithmetic integration step the result needs to be enclosed in a box for the subsequent integration step. When iterated, this causes the `wrapping effect': the final enclosure is much larger than the actual image of the initial set by the flow. \label{fig:wrapping}}
\end{figure}
{\bf Controlling the wrapping effect.} We integrate the ODEs using interval arithmetic. We recall that in interval arithmetic numbers are represented by intervals and vectors by `boxes' which are cartesian products of intervals. When performing arithmetic operations on intervals, an interval enclosure is created by computing the minimum and maximum possible values that can arise from the operation. Repeated algebraic operations can lead to blowup of the interval enclosures. One of the main sources of blowup of the interval enclosures of solutions  of ODEs is due to the wrapping effect.
%
%When integrating an ODE in interval arithmetic one of the main sources of blowup of the interval enclosures of solutions is due to the wrapping effect.
(See Figure \ref{fig:wrapping}.) We overcome this effect by carefully choosing our local coordinates. (See Figure \ref{fig:wrapping-control}.) Our local coordinates $\left( u,s,\alpha,I\right) $ are chosen in such a way that we can exploit the hyperbolicity on coordinates $(u,s)$ for the construction of correctly aligned windows. This allows us to `trim' the sets along the `unstable' coordinate $u$. This way we avoid the wrapping effect on $(u,s)$ since we can obtain trajectories passing through the required domains by using topological shadowing arguments. The coordinate $I$ is almost preserved, so the wrapping effect along it is minimal. We thus essentially reduce the wrapping effect to the coordinate $\alpha $, on which it is unavoidable due to the `twist' dynamics in $\left( \alpha ,I\right)$ coordinates.

{\bf Dealing with center coordinates. }
The method of correctly aligned windows enables the shadowing of sequences of windows of arbitrary length in a hyperbolic-like setting \cite{GideaZ04a, zgliczy2009covering}, where some coordinates are topologically expanding while the others are topologically contracting.  
%The method of correctly aligned windows allows for the shadowing of sequences of windows of arbitrary length in a `hyperbolic' setting  \cite{GideaZ04a, zgliczy2009covering}, where part of the coordinates are expanding and the remaining are contracting. 
 In our system, however, center coordinates $\alpha, I$ are also present. To address this, we use the enhanced method from Theorem \ref{th:diffusion-mechanism} which includes an appropriate overlap of the windows, cone conditions and the condition that sequences of windows return to the same strip. As a result, we can shadow orbits of arbitrary length by verifying only a finite number of conditions. \smallskip

\end{remark}
\begin{remark}
We have chosen the Neptune-Triton system as the example of application of our method, since it exhibits the smallest eccentricity in our solar system. The method is not restricted though to this particular system and can be applied for the three body problem with other mass parameters.
\end{remark}

\section{Acknowledgements}
We thank Jean-Pierre Marco for suggesting to work on this problem.

M.C. was partially supported by the NCN grant 2021/41/B/ST1/00407.  M.G. was partially supported by the NSF grant  DMS-2307718.

The work was conducted as part of the American Institute of Mathematics SQuaREs program ``Arnold diffusion for non-convex Hamiltonians''.

We gratefully acknowledge Polish high-performance computing infrastructure PLGrid (HPC Centers: ACK Cyfronet AGH) for providing computer facilities and support within computational grant no. PLG/2024/017471. 
%TCIDATA{Version=5.00.0.2606}
%TCIDATA{LaTeXparent=0,0,MMFedit.tex}

\section{Appendix}\label{sec:proof-K-reduced}

The points $\{q_i\}_{i=0,\ldots,122}$ used for our computation are written out in Table \ref{table:1}. We write out only some of the points and focus on the coordinates $\bar r_2, \bar \phi_2, \bar R_2, \bar \Phi_2$, since $\pi_{\alpha, I}(q_i) = (0, \frac{1}{10})$. The rest of the points which are not written out in the tables can be recovered from the symmetry 
\[S(\alpha,I,\bar r_2, \bar \phi_2, \bar R_2, \bar \Phi_2) = (\alpha,I,\bar r_2, 2 \pi - \bar \phi_2, -\bar R_2, \bar \Phi_2) \]
using
\begin{align*}
q_{49+i} & = S(q_{49-i})  &\qquad \mbox{for }& i=1,\ldots, 49, \\
q_{104+i} & = S(q_{104-i})  & \mbox{for }& i=1,\ldots, 6, \\
q_{110+i} & = q_{98+i}  & \mbox{for }& i=1,\ldots, 12.\\
\end{align*}

\begin{table}
{\scriptsize
\begin{center}
\begin{tabular}{l l l l l} 
\hline 
$i$	& $\bar r_2$ &	$\bar \phi_2$	& $\bar R_2$	& 	$\bar \Phi_2$ \\
\hline
\hline
0	& 0.951431 &	$\pi$	& 0.0	& 	0.9708302420604081 \\
1	& 0.9539204480706741	& 	3.159254701069249	& 	0.0186158135736214	& 	0.9668474615660276 \\
2	& 0.9601873124783044	& 	3.172250096823662	& 	0.02852852583819738	& 	0.9564332083419783 \\
3	& 0.967493801812415	& 	3.177846444678873	& 	0.02670528162497521	& 	0.9416472041564166 \\
4	& 0.9729596393554711	& 	3.175296379680509	& 	0.01465092366257321	& 	0.9222607238842714 \\
5	& 0.9745898039665974	& 	3.164728176455792	& 	-0.001591808033516218	& 	0.895706635519785 \\
6	& 0.9733881191373061	& 	3.146435648989304	& 	-0.003169800812491891	& 	0.8686875637537721 \\
7	& 0.9741238069976927	& 	3.125826557554216	& 	0.005765322398740717	& 	0.8817054372642115 \\
8	& 0.9742209302457691	& 	3.111664191021475	& 	-0.007155957014991104	& 	0.9113472053549975 \\
9	& 0.9703366991844445	& 	3.105482131181055	& 	-0.02233203643764289	& 	0.9335983435902521 \\
0	& 0.9635319628039253	& 	3.107417293907272	& 	-0.02914538215731054	& 	0.9502612319902282 \\
11	& 0.9564226384989571	& 	3.117279873365691	& 	-0.02448029755015469	& 	0.962793200457194 \\
12	& 0.9519546686579373	& 	3.133325942104626	& 	-0.009010018523418055	& 	0.9699935821371525 \\
13	& 0.9522063234827943	& 	3.151626184587504	& 	0.01088992228421495	& 	0.9695915101384193 \\
14	& 0.9570542926235974	& 	3.167201689132748	& 	0.02547157599021309	& 	0.9617535015908595 \\
15	& 0.9642666648638679	& 	3.176323954650804	& 	0.02895837980093849	& 	0.9487998419503282 \\
16	& 0.9708861272689715	& 	3.177449829373493	& 	0.02111957672906738	& 	0.9316933124456364 \\
17	& 0.9743805622245093	& 	3.170476426834265	& 	0.00547455730040294	& 	0.9087226919803745 \\
18	& 0.9739745525817392	& 	3.155532243387512	& 	-0.005992738978376703	& 	0.8787688163721209 \\
19	& 0.9734829664981886	& 	3.134664448438742	& 	0.004294857875668697	& 	0.8702241290295786 \\
20	& 0.9746125652971285	& 	3.117055971702828	& 	0.0001882116925898782	& 	0.8987011663748681 \\
21	& 0.9725697307474053	& 	3.107276640943593	& 	-0.01618687809052623	& 	0.9244462493696631 \\
22	& 0.9668098042774054	& 	3.10553038927528	& 	-0.02739074814395064	& 	0.9432890195519533 \\
23	& 0.9594717176625115	& 	3.111917295098163	& 	-0.02806594457228648	& 	0.9576780564117845 \\
24	& 0.9534693760344695	& 	3.125534052495292	& 	-0.01705902275863833	& 	0.9675712131571056 \\
25	& 0.9514565600444405	& 	3.14342409864074	& 	0.002012892348322641	& 	0.9707896847490411 \\
26	& 0.9544099880755311	& 	3.160810595242435	& 	0.02008201171235922	& 	0.9660607239147455 \\
27	& 0.9609134515203661	& 	3.173157900058608	& 	0.02887335983514226	& 	0.9551477776528213 \\
28	& 0.9681596529645558	& 	3.177957108944551	& 	0.02591782843624878	& 	0.939962325792649 \\
29	& 0.9733107138550087	& 	3.174605092129972	& 	0.01306440560581605	& 	0.9200026559471112 \\
30	& 0.9745341918710951	& 	3.163246581886198	& 	-0.002857470435922413	& 	0.8926510112783124 \\
31	& 0.9733268022753704	& 	3.144316300522308	& 	-0.001836487978919399	& 	0.8676638511669957 \\
32	& 0.9742707050785351	& 	3.124068033360436	& 	0.005275578216873134	& 	0.8847596058344672 \\
33	& 0.97403248639194	& 	3.110693122791548	& 	-0.008810627726748873	& 	0.9139151161252327 \\
34	& 0.9697807491022538	& 	3.105304654383899	& 	-0.02341467351562676	& 	0.9354943638726507 \\
35	& 0.9628309819242602	& 	3.108047002226755	& 	-0.02914976920130642	& 	0.9517475802058978 \\
36	& 0.9558771097285511	& 	3.118628450070806	& 	-0.0232633191963732		& 0.9638736249325024 \\
37	& 0.951851067533996	& 	3.135091912666698	& 	-0.006872125491257036	& 	0.9704440980903218 \\
38	& 0.9526783722356278	& 	3.153301901542159	& 	0.01308931322035655		& 0.9692439624740564 \\
39	& 0.9580045323843621	& 	3.168265996045714	& 	0.02692277350826037		& 0.9607036048213 \\
40	& 0.9654588549794331	& 	3.176406583534748	& 	0.02939585829553249		& 0.9470695477467663 \\
41	& 0.9720733280051895	& 	3.176288584081977	& 	0.02070738011597229		& 0.9288584052729757 \\
42	& 0.9754489288675509	& 	3.167695945053513	& 	0.005362580545702418	& 	0.9034591243030751 \\
43	& 0.9757198353838055	& 	3.150229674858085	& 	0.002089191280482759	& 	0.8704805851302668 \\
44	& 0.979997693963694	& 	3.127362295688045	& 	0.03174880096684519	& 	0.8829567796566783 \\
45	& 0.9886836885022908	& 	3.112066220553185	& 	0.03139149780771592	& 	0.9384147476817041 \\
46	& 0.995569802834356	& 	3.106072624721055	& 	0.02371443990230434	& 	0.9833045060744217 \\
47	& 1.001781832624171	& 	3.107875369711822	& 	0.02702658026245029	& 	1.026089464871101 \\
48	& 1.010226117414668	& 	3.118309744439625	& 	0.03894522167051586	& 	1.08427674685978 \\
49	& 1.017658989135204	& 	$\pi$	& 	0.0	& 	1.152515980877304 \\
\hline
\end{tabular}
\end{center}}
{\scriptsize
\begin{center}
\begin{tabular}{l l l l l} 
\hline
$i$	& $\bar r_2$ &	$\bar \phi_2$	& $\bar R_2$	& 	$\bar \Phi_2$ \\
\hline
\hline
98	& 0.951431 &	$\pi$	& 0.0	& 	0.9708302420604081 \\
99	& 0.9541018558759905	& 	3.159851225733783	& 	0.01918332721652144	& 	0.9665560278968106 \\
100	& 0.9607376651678461	& 	3.172946398212285	& 	0.02880074318500662	& 	0.9554608638692007 \\
101	& 0.9682500265643872	& 	3.177966658715111	& 	0.02579970281769128	& 	0.9397251579000065 \\
102	& 0.9734758397477414	& 	3.174217348500773	& 	0.01222315099977105	& 	0.9188040423641367 \\
103	& 0.97445668513648	& 	3.161852531863615	& 	-0.003855752742117606	& 	0.8898789259180572 \\
104	& 0.9732942397609048	& 	$\pi$	 	& 0.0	& 	0.8671880572332963 \\
\hline
\end{tabular}
\end{center}}
\caption{\label{table:1} The points $q_i$ along the homoclinic. }
\end{table}

\begin{proof}[Proof of Lemma \protect\ref{lem:K-reduced}]
We will show that%
\begin{equation}
\bar{K}_{\varepsilon }=\frac{1}{2}\left( \bar{R}_{1}^{2}+\mathcal{C}%
_{\varepsilon }-2\mathcal{A}_{\varepsilon }\right) ,
\label{eq:K-bar-epsilon-temp}
\end{equation}%
for%
\begin{align*}
\mathcal{A}_{\varepsilon }& =\frac{\bar{r}_{1}^{2}}{1+\varepsilon \bar{r}_{1}%
}, \\
\mathcal{B}_{\varepsilon }& =2\mathcal{A}_{\varepsilon }+\left( \varepsilon 
\mathcal{A}_{\varepsilon }-\bar{r}_{1}\right) ^{2}, \\
\mathcal{C}_{\varepsilon }& =\left( \varepsilon \bar{\Phi}_{2}\right) ^{2}-2%
\bar{r}_{1}\varepsilon \bar{\Phi}_{2}\left( \varepsilon ^{2}\bar{\Phi}%
_{2}-2\right) +\mathcal{B}_{\varepsilon }\left( 1+\varepsilon ^{2}\bar{\Phi}%
_{2}\left( \varepsilon ^{2}\bar{\Phi}_{2}-2\right) \right) .
\end{align*}%
We compute%
\begin{multline*}
\frac{1}{r_{1}}=\frac{1}{1+\varepsilon \bar{r}_{1}}=1-\varepsilon \bar{r}%
_{1}+\left( \frac{1}{1+\varepsilon \bar{r}_{1}}-\left( 1-\varepsilon \bar{r}%
_{1}\right) \right)  \\
=1-\varepsilon \bar{r}_{1}+\frac{\varepsilon ^{2}\bar{r}_{1}^{2}}{%
1+\varepsilon \bar{r}_{1}}=1-\varepsilon \bar{r}_{1}+\varepsilon ^{2}%
\mathcal{A}_{\varepsilon },
\end{multline*}%
which gives us%
\begin{align*}
\frac{1}{r_{1}^{2}}& =\left( 1-\varepsilon \bar{r}_{1}+\varepsilon ^{2}%
\mathcal{A}_{\varepsilon }\right) ^{2} \\
& =1+2\left( -\varepsilon \bar{r}_{1}+\varepsilon ^{2}\mathcal{A}%
_{\varepsilon }\right) +\left( -\varepsilon \bar{r}_{1}+\varepsilon ^{2}%
\mathcal{A}_{\varepsilon }\right) ^{2} \\
& =1-2\varepsilon \bar{r}_{1}+\varepsilon ^{2}\mathcal{B}_{\varepsilon }.
\end{align*}%
We now have%
\begin{align*}
\left( \frac{1-\Phi _{2}}{r_{1}}\right) ^{2}& =\left( 1-2\varepsilon \bar{r}%
_{1}+\varepsilon ^{2}\mathcal{B}_{\varepsilon }\right) \left( 1-\varepsilon
^{2}\bar{\Phi}_{2}\right) ^{2} \\
& =\left( 1-2\varepsilon \bar{r}_{1}+\varepsilon ^{2}\mathcal{B}%
_{\varepsilon }\right) \left( 1+\varepsilon ^{2}\bar{\Phi}_{2}\left(
\varepsilon ^{2}\bar{\Phi}_{2}-2\right) \right)  \\
& =1-2\varepsilon \bar{r}_{1}-2\varepsilon ^{2}\bar{\Phi}_{2}+\varepsilon
^{2}\mathcal{C}_{\varepsilon },
\end{align*}%
so we see that%
\begin{align*}
K=& \frac{1}{2}\left( \varepsilon ^{2}\bar{R}_{1}^{2}+1+2\varepsilon \left( 
\bar{\Phi}_{1}-\bar{r}_{1}\right) -\varepsilon ^{2}2\bar{\Phi}%
_{2}+\varepsilon ^{2}\mathcal{C}_{\varepsilon }\right) -\left( 1+\varepsilon 
\bar{\Phi}_{1}-\varepsilon ^{2}\bar{\Phi}_{2}\right)  \\
& -\left( 1-\varepsilon \bar{r}_{1}+\varepsilon ^{2}\mathcal{A}\right)  \\
=& \varepsilon ^{2}\bar{K}_{\varepsilon }-\frac{3}{2},
\end{align*}%
as required. The fact that (\ref{eq:K-bar-epsilon-temp}) simplifies to (\ref%
{eq:K-bar-epsilon}) follows from an elementary, though slightly laborious
computation.
\end{proof}

\bibliographystyle{unsrt}
\bibliography{papers,diffusion}

\end{document}